\documentclass[12pt]{article}

\def\calx{{\mathcal{X}}}
\def\caly{{\mathcal{Y}}}

\def\({\left(}
\def\){\right)}

\def\bk{\bigskip }

\def\sk{\smallskip }

\def\n{\noindent }
\def\dd{\displaystyle}

\def\D{{\Delta}}
\def\p{\partial}

\def\barr{\begin{array}}
\def\earr{\end{array}}

\def\bit{\begin{itemize}}

\def\eit{\end{itemize}}

\def\D{{\Delta}}

\usepackage{amsmath, amsthm, amscd, amsfonts, amssymb}
\usepackage[usenames]{color}

\usepackage{mathrsfs}

\numberwithin{equation}{section} 

\newtheorem{theorem}{Theorem}[section]

\newtheorem{corollary}[theorem]{Corollary}

\newtheorem{lemma}[theorem]{Lemma}
\theoremstyle{definition}

\newtheorem{example}[theorem]{Example}

\newtheorem{remark}[theorem]{Remark}

\def\S{Schr\"o\-din\-ger}

\def\calx{\mathcal{X}}

\def\bbe{{\mathbb{E}}}
\def\bbr{{\mathbb{R}}}
\def\bbp{{\mathbb{P}}}

\def\1{^{-1}}

\def\calx{{\mathcal{X}}}

\def\9{{\infty}}

\def\lbb{{\lambda}}
\def\a{{\alpha}}
\def\b{{\beta}}
\def\na{{\nabla}}
\def\g{{\gamma}}
\def\wt{\widetilde}

\def\vf{{\varphi}}

\def\pt{{\partial_t}}
\def\D{{\Delta}}

\def\barr{\begin{array}}
\def\earr{\end{array}}

\def\dd{\displaystyle}
\def\bk{\bigskip }
\def\sk{\smallskip}
\def\n{\noindent }

\def\({\left(}
\def\){\right)}
\def\<{\left<}
\def\>{\right>}

\def\wt{\widetilde}
\def\wh{\widehat}
\def\ol{\overline}
\def\ve{{\varepsilon}}

\def\pxi{{\partial_\xi}}
\def\px{{\partial_x}}

\begin{document}

\begin{center}
{\Large{\bf Strichartz and local smoothing estimates for stochastic dispersive equations with linear multiplicative noise}}
\bigskip\bk

{{\large{\bf Deng Zhang}}\footnote{
School of mathematical sciences,
Shanghai Jiao Tong University, 200240 Shanghai, China. E-mail address: dzhang@sjtu.edu.cn}}
\end{center}

\bk\bk\bk

\begin{quote}
\n{\small{\bf Abstract.}
We study a quite general class of stochastic dispersive equations with linear
multiplicative noise,
including especially the Schr\"odinger and Airy equations.
The pathwise Strichartz and local smoothing estimates
are derived here in both the conservative and non-conservative case.
In particular, we obtain the $\bbp$-integrability of constants
in these estimates,
where $\bbp$ is the underlying probability measure.
Several applications are given to nonlinear problems,
including  local well-posedness
of stochastic nonlinear  Schr\"odinger equations with
variable coefficients and lower order perturbations,
integrability of global solutions
to stochastic nonlinear Schr\"odinger equations with constant coefficients.
As another consequence, we prove
as well the large deviation principle for the small noise asymptotics.  }


\end{quote}

\vfill

\section{Introduction and main results}  \label{Sec-Intro}

We are concerned with the stochastic dispersive equation with linear multiplicative noise
\begin{align} \label{equa-dis}
      &dX(t)= i P(x,D) X(t)dt  + F(t) dt  - \mu X(t)dt +X(t)dW,\ \ t\in (0,T),  \\
      &X(0)=X_0.  \nonumber
\end{align}
Here,  $X$ is a complex-valued function  on $[0,T] \times \mathbb{R}^d$, $T\in (0,\9)$,
$P(x,D)$ is a pseudo-differential operator  of order $m$, $m\geq 2$,
$D_j = - i\p_{x_j}$, $D= (D_1, \cdots, D_d)$.
The term $W$ is a colored Wiener process of the form
\begin{equation} \label{W}
     W(t,x)=\sum^N\limits_{j=1}\mu_j e_j(x)\beta_j(t),\ t\geq
     0,\ x\in \mathbb{R}^d,
\end{equation}
and $\mu(x)=\frac{1}{2}\sum^N_{j=1}|\mu_j|^2e_j^2(x)$,
where
$\mu_j \in \mathbb{C}$,
$e_j$ are real-valued functions,
and $\beta_j(t)$ are independent real Brownian motions on a probability space
$(\Omega,\mathcal{F},\mathbb{P})$ with natural filtration
$(\mathcal{F}_t)_{t\geq 0}$.
For simplicity, we consider $N<\9$,
but the arguments in this paper extend also to the case where $N = \9$.
We assume that $X_0$ is $\mathcal{F}_0$ measurable and
$F$ is $\{\mathcal{F}_t\}$-adapted throughout this paper.

Stochastic dispersive equations arise in various fields of physics.
An important model is the stochastic nonlinear Schr\"odinger equation where
$P(x,D) = - \D$, $F= - \lbb i |X|^{\a-1} X$, $\lbb =\pm 1$ and $\a\in (1,\9)$.
This equation is proposed as a model for the propagation of nonlinear dispersive waves in nonlinear or
random media,
the coefficient $\lbb =1$ (resp. $-1$) corresponds to the focusing (resp. defocusing) case.
See e.g. \cite{BCIRG94,BD99,BD03}.
In particular, when $Re \mu_j=0$, $1\leq j\leq N$,
$-\mu X dt +  X dW$
is indeed the Stratonovitch product $ X \circ d W$,
and the mass of the homogeneous solution is {\it pathwisely} conserved.
This case will be called the {\it conservative case} in this paper.
In the {\it non-conservative case},
i.e. ${\rm Re}\mu_j \not = 0$ for some $1\leq j\leq N$,
this equation plays an important role in the application to open quantum systems.
See e.g. \cite{BG09,BH95,H96}.
An important feature in this case is, that
the mass of the homogeneous solution is
no longer a constant,
but a positive continuous martingale,
which implies conservation in mean norm square
which is crucial to define the so called ``physical probability law'' (see \cite{BG09}).
For other dispersive type equations,
see e.g. \cite{BD07} for  the
stochastic Korteweg-de Vries (KdV) equation where
$P(x,D)=D^3$ and $F= \frac 12 \p_xX^2$,
and \cite{CG14, CG15, BD10, DT11} for Schr\"odinger and KdV equations
with modulated dispersion.

Unlike the usual parabolic case,
the principle operator of a dispersive equation usually generates a unitary group in the
standard $L^2$ space.
Thus, a global smoothing effect is excluded
in Sobolev spaces $H^s(\bbr^d)$, $s>0$,
which is the source of many difficulties to study nonlinear problems.
Furthermore,
although the principle operator is monotone,
the variational approach (see e.g. \cite{LR15})
is not applicable to stochastic dispersive equations,
due to the lack of coercivity of the principle operator.

Here we shall study Strichartz and local smoothing estimates
for stochastic dispersion equations,
which are two most stable ways of measuring dispersion
and play an important role to study nonlinear problems.

The Strichartz estimates give  space-time integrability of solutions,
while the local smoothing estimates allow to gain $(m-1)/2$ derivatives of solutions
on every bounded domains.
We refer to \cite{KT98,KPV91, MMT08,RZ05,RS04,T08} for Strichartz estimates
and \cite{C02,D94,D96,D00,K13,KPV91,KPV04,MMT08} for local smoothing estimates in the deterministic case.

For stochastic nonlinear Schr\"odinger equations,
the global well-posedness was first studied in \cite{BD99,BD03}
for general multiplicative noises.
The key Strichartz estimates for the stochastic convolution
were proved there by using the theory of $\g$-radonifying operators,
which, as the role of Hilbert-Schmidt operators on Hilbert spaces,
allows to treat noises in Banach spaces.
An improved stochastic Strichartz estimates was
proved in \cite{BM14}, based on which global well-posedness
was obtained
on a two-dimensional compact manifold.
See  \cite{H16} for the global well-posedness in the full mass subcritical case
via the stochastic Strichartz estimates in \cite{BM14}.
See also the recent work \cite{BHW17} for martingale solutions
in the energy space on compact manifolds.

Recently, using a different approach based on the rescaling transformation (see \eqref{rescal} below),
global well-posedness for stochastic nonlinear Schr\"odinger equations with linear multiplicative noise
has been proved in the optimal mass and energy subcritical cases (\cite{BRZ14,BRZ16}),
where the key role is played  by the pathwise Strichartz and local smoothing estimates.
It should also be mentioned that,
the rescaling approach is quite robust and fits well with
the theory of maximal operators.
In particular,
solutions continuously depend on the initial condition {\it pathwisely}
and satisfy the strict cocycle property,
so give rise to stochastic dynamical systems
(see \cite{A98}).
We refer to \cite{BRZ16.1} for stochastic logarithmic Schr\"odinger equations,
\cite{BRZ14.2} for noise effect in the non-conservative case,
and \cite{BRZ16.2} for optimal control problems.

In this paper, we prove the  pathwise  Strichartz
and local smoothing estimates
for quite general stochastic dispersive equations with linear multiplicative noise
in a uniform manner,
including especially the Schr\"odinger and Airy equations.

Moreover, motivated by scattering and optimal control problems,
we also obtain explicit  upper bounds
and $\bbp$-integrability  of constants in
these estimates.
In particular, for the homogeneous Schr\"odinger and Airy equations,
the constants in the local smoothing estimates are exponentially $\bbp$-integrable in the conservative case.

Several applications are given to nonlinear problems.
Pathwise local well-posedness is proved
for stochastic nonlinear Schr\"odinger equations
with  variable coefficients and also with lower order perturbations.
Moreover,  the $\bbp$-integrability of global solutions
are obtained for the stochastic nonlinear Schr\"odinger equation
mentioned above
in both the mass and energy subcritical cases,
which can  be viewed as a complement to \cite{BRZ14,BRZ16, BRZ14.2}
and is of importance for optimal control problems (see \cite{BRZ16.2}).

Another application we obtain in this paper is that
the large deviation principle for the small noise asymptotics
for general linear stochastic dispersive equations,
as well as stochastic nonlinear Schr\"odinger equations with variable coefficients,
in the conservative case. \\

{\bf Notations.}
For any $x=(x_1,\cdots,x_d) \in \bbr^d$
and any multi-index $\a=(\a_1,\cdots, \a_d)$,
$\<x\>=(1+|x|^2)^{1/2}$,
$|\a|= \sum_{j=1}^d \a_j$,
$\partial_x^\a=\partial_{x_1}^{\a_1}\cdots \partial_{x_d}^{\a_d}$,
$\<\partial_x\>=(I-\Delta_x)^{1/2}$.
Let $D_{x_j} = -i \partial_{x_j}$,
$D_x^\a$ is defined similarly.
We will use the notation $\xi$ for the phase variable.

Given $1\leq p \leq \9$, $p'$ is the conjugate number, i.e., $1/{p'}+1/p=1$.
$L^p = L^p(\bbr^d)$ is the space of $p$-integrable complex functions with the norm $|\cdot|_{L^p}$.
In particular, $L^2$ is the Hilbert space with the inner product $\<u,v\>_2 = \int u(x)\ol{v(x)}dx$,
and $|\cdot|_2 = |\cdot|_{L^2}$.  As usual,
$\dot{W}^{s,p}= |D|^{-s}L^p(\bbr^d)$,
$W^{s,p}= \<D\>^{-s}L^p(\bbr^d)$,
and $H^1 = W^{1,2}$.
Let $\mathscr{S}$ denote the space of rapid decreasing functions
and $\mathscr{S}'$ the dual space of $\mathscr{S}$.
For any $f\in \mathscr{S}$, $\wh{f}$
is the Fourier transform of $f$,
i.e. $\wh{f}(\xi) = \int e^{-ix\cdot \xi} f(x)dx$.

For any Banach space $\calx$,
$L^p(0,T; \calx)$ is the space of $p$-integrable $\calx$-valued functions
with the norm $\|\cdot\|_{L^p(0,T; \calx)}$,
and $C([0,T];\calx)$ is the space of continuous
$\calx$-valued functions with the super norm in $t$.
For two Banach spaces $\calx, \caly$,
$\mathcal{L}(\calx,\caly)$ is the space of  linear continuous operators from $\calx$ to $\caly$,
and $\mathcal{L}(\calx)= \mathcal{L}(\calx,\calx)$.

Throughout this paper, we use $C(\cdots)$ for various constants that may
change from line to line.

\section{Formulations of main results} \label{Sec-Main}

To begin with, let us introduce suitable spaces for Strichartz and local smoothing estimates.

Set
$H^{s}_{\delta}:= \{v\in \mathscr{S}': \<x\>^{\delta}(I-\Delta)^{s/2} v \in L^2\}$
with the norm
$\|v\|_{H^{s}_{\delta}}= \|\<x\>^{\delta}(I-\Delta)^{s/2}v\|_{2}$.
For any $s\in [0,\9)$, $p,q\in [2,\9]$,
set
$$\calx_{T,s,p,q} := L^q(0,T;\dot{W}^{s,p}) \cap L^2(0,T;H^{\frac{m-1}{2}}_{-1}),$$
equipped with the norm
$\|u\|_{\calx_{T,s,p,q}} = \|u\|_{L^q(0,T; \dot{W}^{s,p})} + \|u\|_{L^2(0,T; H^{\frac{m-1}{2}}_{-1})}$.

The dual space of $\calx_{T,s,p,q}$ is denoted by
$\calx'_{T,-s,p',q'} := L^{q'}(0,T;\dot{W}^{-s,p'}) + L^2(0,T;H^{-\frac{m-1}{2}}_{1})$,
and $\|u\|_{\calx'_{T,-s,p',q'}} = \inf \{\|u_1\|_{L^{p'}(0,T; \dot{W}^{-s,p'})}+ \|u_2\|_{L^2(0,T; H^{-(m-1)/2}_{1})}: u=u_1+u_2,
u_1 \in L^{p'}(0,T; \dot{W}^{-s,p'}), u_2 \in L^2(0,T; H^{-(m-1)/2}_{1})\}$.

We say that
$a\in C^\9(\bbr^d \times \bbr^d)$ is a symbol of class $S^m$,
if for any multi-indices $\a,\beta \in \mathbb{N}^d$,
$|\partial^\a_\xi\partial^\beta_x a(x,\xi)| \leq C_{\a,\beta} \<\xi\>^{m-|\a|}$.
The semi-norms $|a|_{S^m}^{(l)}$, $l\in \mathbb{N}$, are defined by
$$|a|_{S^m}^{(l)} = \max_{|\a+\beta|\leq l} \sup_{\bbr^{2d}}
                        \{ |\partial^\a_\xi \partial^\beta_x a(x,\xi)|\<\xi\>^{-(m-|\a|)}\}.$$
Moreover, $\Psi_a$ (or $a(x,D)$) denotes the pseudo-differential operator
associated with the symbol $a(x,\xi)$, i.e.,
\begin{align*}
   \Psi_a v(x) = a(x,D) v(x) = (2\pi)^{-d} \int e^{ix \cdot \xi} a(x,\xi) \wh{v}(\xi) d\xi,\ \ v\in \mathscr{S}.
\end{align*}
In this case, we write $\Psi_a \in S^m$ (or $a(x,D) \in S^m$) when no confusion arises.

We shall assume that the principle symbol and the spatial functions in the noise satisfy
the following assumptions:
\begin{enumerate}
  \item[(A0)] $P(x,\xi) = P_1(x,\xi) + i P_2(x,\xi)$, $P_1\in S^m$,
  $P_2\in S^{m-1}$, $P_j(x,\xi)$ are real polynomials of $\xi$, and $P(x,D)$ is self-adjoint.
  \item[(A1)] {\it Asymptotical flatness.}
  For any multi-indices $\a, \beta$,
  \begin{align*}
    |\p_\xi^\a \p_x^\beta P_1(x,\xi)| \leq& C_{\a\beta} \<x\>^{-2} \<\xi\>^{m-|\a|},\ \beta\not=0. \\
    |\p_\xi^\a \p_x^\beta P_2(x,\xi)| \leq& C_{\a\beta} \<x\>^{-2} \<\xi\>^{m-1-|\a|}.
  \end{align*}
  Moreover, for $1\leq j \leq N$, $e_j\in C^\9(\bbr^d)$, and
  \begin{align} \label{e-decay}
      |\p_x^\a  e_j(x)| \leq C_\a \<x\>^{-2},\ \ \a\not =0.
  \end{align}
\item[(A2)] {\it Non-trapping condition.}
  The bicharacteristic flow associated with the principle symbol $P(x,\xi)$ is non-trapped.
  More precisely, let
  $(X,\Xi)(t,x,\xi)$ be a flow generated by the Hamiltonian vector field
  \begin{align*}
    H_P = \sum \p_{\xi_j} P(x,\xi) \p_{x_j} - \p_{x_j} P(x,\xi) \p_{\xi_j} ,
  \end{align*}
  that is, $(X,\Xi)$ is a solution to
  \begin{align*}
    &\frac{dX}{dt} = \na_\xi P(X,\Xi), \ \ X(0)=x,\\
    & \frac{d \Xi}{dt} = (-1) \na_x P(X,\Xi),  \ \ \Xi(0)=\xi.
  \end{align*}
  Then, for any $(x,\xi) \in T^*\bbr^{d} / \{0\}$, $|X(t)| \to \9$,
  as $t\to \pm \9$.
  \item[(A3)] {\it Strichartz-type estimate.}
  There exists $(s,p,q)\in[0,\9) \times [2,\9) \times (2,\9]$ such that
  \begin{align} \label{Stri-P}
      \|e^{i\cdot P} u \|_{L^{q}(0,\9; \dot{W}^{s,p})}
      \leq C |u|_2.
  \end{align}
  The triple $(s,p,q)$ such that \eqref{Stri-P} holds is called admissible,
  and the set of all admissible triples is denoted by $\mathcal{A}$.
  In particular, $(0,2,\9) \in \mathcal{A}$.
  (Note that, we do note treat the endpoint case $p=\9$, $q=2$ here.)
\end{enumerate}

\begin{remark} \label{Rem-Stri}
Assumptions $(A0)$--$(A2)$ are mainly required for
local smoothing estimates.
See Theorem \ref{Thm-Loc-Smooth} and Lemma \ref{Lem-Doi} below.
See also \cite{C02} for more general conditions for local smoothing estimates.
We mention that, for $P$ elliptic with variable coefficients
under some appropriate flatness conditions at infinity,
it was proved in \cite{D00}
that the local smoothing estimate is equivalent to
the non-trapping condition of the Hamiltonian flow.
\end{remark}

\begin{remark}
The smoothness of the spatial functions $e_j$
is   assumed for technical reasons, particularly, to perform pseudo-differential calculus.
One can also assume $e_j\in C_b^N(\bbr^d)$ as in \cite{KPV04} with $N$ large enough
such that the pseudo-differential calculus can
be carried out.
One may also weaken the regularity of $e_j$ by using paradifferential calculus as
in \cite{T08,MMT08}.
\end{remark}

\begin{example}
{\it Schr\"odinger operator.} Consider $P(x,D)=\sum_{j,k=1}^d D_j a^{jk}(x) D_k$,
$x\in \bbr^d$, $d\geq 1$,
$a^{jk}$ are real valued, symmetric, and satisfy some appropriate conditions (see Assumptions (B1) and (B2) below).
In this case, $P(x,\xi) = \sum_{j,k=1}^d (a^{jk}(x) \xi_j\xi_k - i \p_ja^{jk}(x) \xi_k)$.
In particular, when $a^{jk} = \delta_{jk}$, $P(x,D) = -\Delta$.
We have the admissible set (see \cite{LP09})
\begin{align} \label{Stri-Sch}
  \mathcal{A}=\bigg\{(s,p,q): s=0, \frac 2q = d(\frac 12 - \frac 1p), (p,q)\in [2,\9)\times(2,\9]  \bigg\}.
\end{align}
The pairs $(p,q)$ in \eqref{Stri-Sch} are called {\it Strichartz pairs} below.
\end{example}

\begin{example}
{\it Airy operator.} Consider $P(x,D)= D^3$, $d=1$,
and so $P(x,\xi) = \xi^3$, $\xi\in \bbr$.
This operator mainly arises in the (generalized) KdV equations, for which we have
(see \cite[(1.5)]{KPV89})
\begin{align*}
   \mathcal{A}
   = \bigg\{ (s,p,q): s= \frac{\theta \a}{2}, p= \frac{2}{1-\theta}, q=\frac{6}{\theta(\a+1)}, (\theta,\a)\in [0,1) \times[0,1/2] \bigg\}.
\end{align*}
\end{example}

\begin{example}
Generalization of the Schr\"odinger operator $P(x,D) = (-\Delta)^{m}$,
$m\in \mathbb{N}$, $m\geq 1$, $d\geq 2$.
We have (see Remark (a) on p.49 in \cite{KPV91})
\begin{align*}
   \mathcal{A} = \bigg\{ (s,p,q): s= \frac \theta 2 d(m-1), p=\frac{2}{1-\theta}, q=\frac{4}{d\theta}, \theta\in [0,\frac 2 d) \bigg\}.
\end{align*}
Moreover, for the generalization of the Airy operator
$P(x,D) = D^m$, $m\in \mathbb{N}$, $m\geq 3$, $d=1$,
we have (see \cite[Theorem 2.1]{KPV91})
\begin{align*}
   \mathcal{A} = \bigg\{ (s,p,q): s= \frac{(m-2)\theta}{4}, p=\frac{2}{1-\theta}, q=\frac{4}{\theta}, \theta\in [0,1) \bigg\}.
\end{align*}
\end{example}

The main result of this paper is formulated below.
\begin{theorem} \label{Thm-Stri-L2}
Assume $(A0)$-$(A3)$. Let  $(s_1,p_1,q_1),(s_2,p_2,q_2) \in \mathcal{A}$
be any admissible triples.
Then, we have

$(i)$. For any $\mathcal{F}_0$-measurable $X_0\in L^2$,
$\{\mathcal{F}_t\}$-adapted $F\in \calx'_{T,-s_2,p_2',q_2'}$, $\bbp$-a.s.,
the solution $X$ to \eqref{equa-dis} satisfies that $\bbp$-a.s.
\begin{align} \label{Stri-X}
   \|e^{-\Phi(W)}X\|_{\calx_{T,s_1,p_1,q_1}} \leq  D(e,e^*, T)  (|X_0|_2+ \|e^{-\Phi(W)}F\|_{\calx'_{T,-s_2,p_2',q_2'}}).
\end{align}
Here,
\begin{align} \label{PhiW}
\Phi(W) = W - \frac 12 t \sum_{j=1}^N (|\mu_j|^2 + \mu_j^2) e^2_j,
\end{align}
and
\begin{align} \label{D}
   D(e,e^*, T)
   =& C_1 (1+T)^{\frac 52} (1+e_T+|e|_{L^1(0,T)})^{\frac 12}
     (1+e^*_T+|e^*|_{L^1(0,T)})^{\frac 12}   \nonumber \\
     &\cdot      (1+\beta^*_T + T)^{\frac 32(m^2-m+4)l} e^{C_2 (1+\beta^*_T+T)}
\end{align}
for some $l\geq 1$,
where $\beta^*_T= \sup_{t\in [0,T]}|\beta(t)|$,
$C_1>0$, $C_2\geq 0$,
and $e$ and $e^*$ are adapted processes as in \eqref{U-l2} and \eqref{U*-l2} below respectively.
In particular, one can take $l=e=e^*=1$,  $C_2=0$ in the conservative case.

$(ii)$.
Assume in addition that $P(x,\xi)=P(\xi)$,
$X_0\in H^1$,
$\p_{x_j} F \in \mathcal{X}_{-s_3,p'_3,q'_3}$
for some $(s_3,p_3,q_3)\in \mathcal{A}$, $1\leq j\leq d$, $\bbp$-a.s.
Then, the solution to \eqref{equa-dis} satisfies  that $\bbp$-a.s.,
for any $1\leq j\leq d$,
\begin{align} \label{Stri-naX}
   \|\p_{x_j} (e^{\Phi(W)}X) \|_{\calx_{T,s_1,p_1,q_1}}
   \leq&(1+\beta^*_T + T)^m D^2(e,e^*, T)
         \bigg(|X_0|_2  + \|e^{-\Phi(W)}F\|_{\calx'_{T,-s_2,p_2',q_2'}}\nonumber  \\
       &    + \|\p_{x_j} (e^{-\Phi(W)}F)\|_{\calx'_{T,-s_3,p_3',q_3'}} \bigg).
\end{align}
\end{theorem}

\begin{remark}
The process $e$ is actually related to the martingale property of homogeneous solutions in the
non-conservative case,
while $e^*$  arises in the duality case and involves
the semi-martingale $M^*$ in \eqref{Mar*} below.
\end{remark}

\begin{remark}
The estimate \eqref{Stri-X} shows that for $\bbp$-a.e. $\omega$,
$X(t,\omega) \in H^{(m-1)/2}_{-1} \cap \dot{W}^{s,p}$  for a.e. $t\in [0,T]$.
Hence, the solution gains $(m-1)/2$ derivatives on any bounded domains
and also gain spatial integrability.
We remark that,
the local smoothing estimate (resp. Strichartz estimate) depends
on the first (resp. second) derivatives of the principle symbol (see \cite{KPV91}).
\end{remark}

\begin{remark} \label{Rem-Integ}
The upper bound in \eqref{Stri-X} (and also \eqref{Stri-naX}) can  be improved in the homogeneous case
(see Theorem \ref{Thm-Loc-Smooth} and Lemma \ref{Homo-U} below),
but we will not seek the optimal upper bound here.
However, we have the $\bbp$-integrability of constants and  solutions
which are important for optimal control problems (see \cite{BRZ16.2}).

More precisely, similarly to \cite[Lemma 3.6]{BRZ16}, we have that for any $1\leq \rho<\9$,
$\bbe  \sup_{t\in[0,T]}$ $(M(t)+M^*(t)) ^\rho \leq C(\rho)$,
where $M$ and $M^*$ are as in \eqref{Mar} and \eqref{Mar*} below respectively.
Thus, it follows from \eqref{D} that
\begin{align*}
   D(e,e^*,T)  \in \bigcap\limits_{1\leq \rho<\9} L^\rho(\Omega).
\end{align*}

Moreover, when $X_0 \in L^\rho(\Omega; L^2)$ and
$F \in L^\rho(\Omega; \calx'_{T,-s_0,p'_0,q'_0})$ for some admissible triple $(s_0,p_0,q_0)\in \mathcal{A}$
and for all $1\leq \rho<\9$,
we have that
$$X \in  \bigcap\limits_{1\leq \rho<\9} \bigcap\limits_{(s,p,q)\in\mathcal{A}} L^\rho(\Omega; \calx_{T,s,p,q}). $$
In particular, by the mild reformulation of \eqref{equa-dis}, this implies
the $\bbp$-integrability of the stochastic convolution, i.e.,
\begin{align*}
   \int_0^\cdot e^{i(\cdot-s)P(x,D)} X(s) dW(s) \in  \bigcap\limits_{1\leq \rho<\9} \bigcap\limits_{(s,p,q)\in\mathcal{A}} L^\rho(\Omega; \calx_{T,s,p,q}).
\end{align*}
\end{remark}

\begin{remark}
The $\bbp$-integrability of solutions to \eqref{equa-dis} can be proved by the
stochastic Strichartz estimate in \cite{BM14}.
Here,
we obtain quantitative information for general dispersive equations with linear multiplicative noise.
As a matter of fact,
we obtain \eqref{Stri-X} in the pathwise way.
In particular, for the Schr\"odinger and Airy operators,
the constants in the homogeneous local smoothing estimates are exponentially $\bbp$-integrable in the conservative case.
See Theorem \ref{Thm-Loc-Smooth} and Remark \ref{Rem-Exp-integ} below.
Moreover,
stochastic dispersive equations with lower order perturbations can be treated here,
and the pathwise estimates are applicable as well to the large deviation principle
for the small noise asymptotics in the conservative case.
See Theorems \ref{Thm-LDP} and  \ref{Thm-LDP-SNLS} below.
\end{remark}

Below we present the Strichartz and local smoothing estimates for stochastic dispersive equations
with lower order perturbations.
\begin{theorem}  \label{Thm-Stri-per}
Consider the equation
\begin{align} \label{equa-per}
   &dX= i P(x,D) Xdt + b(t,x,D)X dt  + F dt  - \mu X(t)dt +X(t)dW,   \\
      &X(0)=X_0.   \nonumber
\end{align}
Here, $P(x,D)$, $\mu$ and $W$ are as in \eqref{equa-dis}.
For each $x,\xi\in \bbr^d$,
$t \mapsto b(t, x,\xi)$ is an adapted continuous process
satisfying that, for any
multi-indices $\a, \beta$,
\begin{align} \label{b}
 \sup_{t\in [0,T]} |\p_\xi^\a \p_x^\beta b(t,x,\xi)| \leq g(T) C_{\a\beta}\<x\>^{-2} \<\xi\>^{m-1-|\a|}, \ \ \bbp-a.s.,
\end{align}
where $t \mapsto g(t)$ is an $\{\mathcal{F}_t\}$-adapted   continuous process.

Then, under Assumptions $(A0)$--$(A3)$,
for any $(s_i, p_i, q_i) \in \mathcal{A}$, $i=1,2$,
and for any $\mathcal{F}_0$-measurable $X_0\in L^2$, $\{\mathcal{F}_t\}$-adapted  $F\in \mathcal{X}'_{T,-s_2,p_2',q_2'}$, $\bbp$-a.s.,
the solution $X$ to \eqref{equa-per} satisfies that $\bbp$-a.s.,
\begin{align}  \label{prop-stri}
   \|e^{-\Phi(W)}X\|_{\mathcal{X}_{T,s_1,p_1,q_1}}
   \leq C_T (|X_0|_2 + \|e^{-\Phi(W)}F\|_{\mathcal{X}'_{T,-s_2,p'_2,q'_2}}),
\end{align}
where $\Phi(W)$ is as in \eqref{PhiW} and $t\mapsto C_t$ is adapted, increasing and  continuous.
\end{theorem}

\begin{remark}
The constant $C_T$ in \eqref{prop-stri} may not be $\bbp$-integrable,
due to the lack of martingale property of homogeneous solutions to \eqref{equa-per} with lower order perturbations in general.
In fact, in the derivation of \eqref{prop-stri}  the Gronwall inequality
will produce a non-integrable double exponential of Brownian motions (see e.g. \eqref{d-exp} below).
\end{remark}

As mentioned above,
several applications are given to nonlinear problems,
which are actually main motivations for the estimates in Theorems \ref{Thm-Stri-L2} and \ref{Thm-Stri-per}.
We  first show
the pathwise local well-posedness for stochastic nonlinear Schr\"odinger equations
with variable coefficients and lower order perturbations,
in the full mass (sub)critical range of the exponents of nonlinearity.
See Theorem \ref{Thm-QSNLS} below.
For the typical stochastic nonlinear Schr\"odinger equation as studied in \cite{BRZ14,BRZ16,BD99,BD03,H16},
we also prove the $\bbp$-integrability of global solutions
which is important for optimal control problems.
See Theorem \ref{Thm-SNLS}.

Moreover,
these estimates apply also to
the large deviation principle for the small noise asymptotics for linear stochastic dispersive equations,
as well as stochastic nonlinear Schr\"odinger equations with variable coefficients,
in the conservative case.
See Theorems \ref{Thm-LDP} and \ref{Thm-LDP-SNLS} below.

The proof of Theorem \ref{Thm-Stri-L2} relies on
the rescaling approach as in \cite{BRZ14,BRZ16,BRZ14.2}.
The rescaling transformation is in fact a Doss-Sussman type transformation in infinite dimension,
which reduces the stochastic dispersive equation \eqref{equa-dis}
to a random equation \eqref{equa-rnls} below with lower order perturbations.

This point of view allows pathwise analysis of stochastic partial differential equations,
including sharp pathwise estimates of stochastic solutions,
path-by-path uniqueness and random attractors.
See e.g. \cite{BDR09,BR13} for stochastic porous media equations
and the total variations flow.
Moreover, the rescaling approach fits quite well with
the theory of maximal monotone operators
and indeed reveals the structure of stochastic equations.
See e.g. \cite{BR15,BRZ16.1}.
We would also like to mention that, the damped effect of the noise in the non-conservative case,
completely different from that in the conservative case,
can be revealed by the rescaling approach (see \cite{BRZ14.2}).

Below we shall see that the rescaling transformation leaves the principle symbol unchanged
but produces lower order perturbations in the resulting random equation.
In the light of this structural feature,
we shall prove Strichartz estimates
via the control of lower order perturbations.

More precisely, let
\begin{align} \label{rescal}
    u := e^{-\Phi(W)}X,
\end{align}
where $\Phi(W)$ is as in \eqref{PhiW}.
By \eqref{equa-dis} we have
\begin{align} \label{equa-rnls}
     \partial_t u (t) = i P_t(x,D) u(t)  + f(t)
\end{align}
with the initial datum $u(0)= X_0$, $ f(t) = e^{-\Phi(W)(t)} F(t)$, and
$$P_t(x,D) = e^{-\Phi(W)(t,x)} P(x,D) e^{\Phi(W)(t,x)}, \ \ t\in[0,T].$$
In particular, $P_0(x,\xi) = P(x,\xi)$.
The equivalence of solutions to  \eqref{equa-dis} and \eqref{equa-rnls}
can be proved similarly as in \cite[Lemma 6.1]{BRZ14}.
Note that
\begin{align} \label{Pt-P}
   P_t = e^{-\Phi(W)(t)} P e^{\Phi(W)(t)}
   = P + e^{-\Phi(W)(t)} [P, e^{\Phi(W)(t)}].
\end{align}
It follows that
\begin{align}\label{equa-rnls2}
   \partial_t u (t) = i P u(t) + ie^{-\Phi(W)(t)} [P, e^{\Phi(W)(t)}] u(t) + f(t),
\end{align}
where $ e^{-\Phi(W)} [P, e^{\Phi(W)}]$
is of lower order $m-1$.
Therefore,
the original  problem is now reduced to this random equation.

\begin{theorem} \label{Thm-Stri-RNLS}
Assume $(A0)$-$(A3)$.
$(i)$. For any $(s_i, p_i,q_i)\in \mathcal{A}$, $i=1,2$,
and for any
$\mathcal{F}_0$-measurable  $u_0\in L^2$ and
$\{\mathcal{F}_t\}$-adapted
$f\in \calx'_{T,-s_2, p_2',q_2'}$, $\bbp$-a.s.,
the solution $u$ to \eqref{equa-rnls} satisfies $\bbp$-a.s. that
\begin{align} \label{l-stri-esti-lps}
  \|u\|_{\calx_{T,s_1,p_1,q_1}}
  \leq  D(e,e^*,T) (|u_0|_2 + \|f\|_{\calx'_{T,-s_2,p_2',q_2'}}),
\end{align}
where $D(e,e^*,T)$ is as in \eqref{D}.

$(ii)$. Assume in addition that $P(x,\xi) = P(\xi)$, $u_0\in H^1$,
and $\partial_{x_j} f \in \mathcal{X}'_{T,-s_3,p_3',q_3'}$ for
some $(s_3,p_3,q_3) \in \mathcal{A}$, $1\leq j\leq d$, $\bbp$-a.s.
Then,
\begin{align} \label{l-stri-esti-wps}
  \|u\|_{\calx_{T,s_1,p_1,q_1}}
  \leq& (1+\beta^*_T + T)^m D(e,e^*,T) \bigg(|u_0|_{H^1} + \|f\|_{\calx'_{T,-s_2,p_2',q_2'}} \nonumber \\
      & + \|\partial_{x_j} f\|_{\calx'_{T,-s_3,p_3',q_3'}} \bigg).
\end{align}
\end{theorem}

\begin{remark}
The proof presented below applies as well to the stochastic dispersive equation of more general form
\begin{align}   \label{equa-gdis}
      dX= i P(x,D) Xdt  + F dt  - \mu Xdt + \lbb X dt + \sum\limits_{j=1}^N F_kXd\beta_k(t),
\end{align}
where $X(0) = X_0$, $\lbb \in \mathbb{C}$,
$\mu= \frac 12 \sum_{j=1}^N |F_k|^2$,
and $F_k$ are adapted complex valued functions on $\bbr^+ \times \bbr^d$.
Under appropriate conditions of $F_k$,
we can perform the rescaling transformation
$u= e^{-\wt{\Phi}(\beta) - \lbb t} X$ with
\begin{align*}
   \wt{\Phi}(\beta)(t,x) = \sum\limits_{j=1}^N \int_0^t F_k(s,x)\beta_k(s)
       - \frac 12 \sum\limits_{j=1}^N \int_0^t (|F_k(s,x)|^2 + F^2_k(s,x) ) ds
\end{align*}
to reduce \eqref{equa-gdis} to the random equation below
\begin{align*}
    \p_t u = i \wt{P}_t(x,D) u + e^{-\wt{\Phi}(\beta)-\lbb t} F,
\end{align*}
where $\wt{P}_t(x,D) = e^{-\wt{\Phi}(\beta)(t)} P(x,D)  e^{\wt{\Phi}(\beta)(t)}$.
Thus, using similar arguments below one can prove Strichartz and local smoothing estimates
for \eqref{equa-gdis}.
\end{remark}

Below we are mainly concerned with the estimates
in Theorems \ref{Thm-Stri-L2}, \ref{Thm-Stri-per}  and \ref{Thm-Stri-RNLS}.
The global well-posedness for \eqref{equa-dis}, \eqref{equa-per} and \eqref{equa-rnls}
can be proved via approximation procedure with smooth initial data and smooth inhomogeneous parts.

There is an extensive literature on
Strichartz estimates for the free group $\{e^{itP(x,D)}\}_{t\in \bbr}$,
of which the standard proof is based on  dispersive estimates,
e.g.,
$|e^{-it\Delta} u_0|_{L^\9} \leq C t^{-d/2} |u_0|_{L^1}$
for the Schr\"odinger operator, and
$|e^{-t\p_x^3} u_0|_{L^\9}$ $\leq C t^{-\frac 13} |u_0|_{L^1}$
for the Airy operator.
See \cite{KT98,LP09,T06}.

However, it is much more difficult to prove
Strichartz estimates for operators with lower order perturbations
and,
as a matter of fact,
dispersive estimates do not hold in general
(see e.g. \cite{T08}).

Inspired by the work \cite{MMT08,RZ05,RS04},
we shall prove Strichartz estimates by using
local smoothing estimates
under appropriate asymptotically flat conditions,
which allow to control lower order perturbations.
For this purpose,
we  first prove the local smoothing estimates
for \eqref{equa-rnls2} in the homogeneous case
(see Theorem \ref{Thm-Loc-Smooth} below),
which actually plays the key role in the proof of Theorem \ref{Thm-Stri-RNLS}.

It should be mentioned that,
local smoothing estimates for more general operators with lower order perturbations were studied in \cite{C02},
where, however, the dependence on time of the constants  is implicit.

In order to obtain the $\bbp$-integrability of constants,
we prove precisely upper bounds of constants,
by estimating the
remainders in the expansion of compositions of pseudo-differential operators
and also the
semi-norms of pseudo-differential operators
(see Lemma \ref{Lem-EstiErr} and Corollary \ref{Cor-EstiErr} below).
We  first treat the easier conservative case in the spirit of \cite{CKS95},
and then, for the non-conservative case,
we perform the energy method to a new equation
as in \cite{KPV04},
combined with the
G{\aa}rding inequality \eqref{garding}  and the interpolation estimate
\eqref{inter} below.
Moreover, the Gronwall inequality used in \cite{C02}
produces a non-integrable double exponential boundedness
of Brownian motions.
Instead, we shall use the martingale property of homogeneous solutions
to obtain the $\bbp$-integrability of constants.

Once  the homogeneous local smoothing estimates are obtained,
by virtue of  Assumption $(A3)$ and the Christ-Kiselev lemma,
we obtain homogeneous Strichartz estimates
and prove Theorem \ref{Thm-Stri-RNLS} by duality arguments,
thereby proving Theorem \ref{Thm-Stri-L2} via the rescaling transformation.

The remainder of this paper is structured as follows. In Section \ref{Sec-Pre}
we review basic results of pseudo-differential operators
and  present  necessary estimates
used in subsequent sections.
Section \ref{Sec-Local-Smooth} is mainly concerned with the homogeneous local smoothing estimates,
and Section \ref{Sec-Stri-Esti} is devoted to the proof of the main results.
In Section \ref{Sec-Appli}, we present several applications concerning stochastic nonlinear
Schr\"odinger equation as well as large deviation principle for small noise asymptotics.
Finally, for simplicity of exposition, some technical computations are postponed to the Appendix, i.e., Section \ref{Sec-App}.

\section{Preliminary} \label{Sec-Pre}

We first review some basic results  of pseudo-differential operators.
For more details we refer to \cite{KPV04,K81,LP09,T00}.

\begin{lemma} \label{Lem-Err} (\cite[Theorem 2.6, Theorem 3.1]{K81})
Let $a_i\in S^{m_i}$, $i=1,2$. Then, $\Psi_{a_1} \circ \Psi_{a_2} = \Psi_a$ with
\begin{align*}
   a(x,\xi)
   = (2\pi)^{-d} \iint e^{-iy\cdot \eta} a_1(x,\xi+\eta) a_2(x+y,\xi) dy d\eta \in S^{m_1+m_2},
\end{align*}
Moreover, we have the expansion
\begin{align*}
  a(x,\xi) = \sum\limits_{|\a|<n}
             \frac{1}{\a!}\partial^\a_\xi a_1(x,\xi) D^\a_x a_2(x,\xi)
              + n\sum_{|\g|=n} \int_0^1 \frac{(1-\theta)^{n-1}}{\g !} r_{\g,\theta}(x,\xi) d\theta,
\end{align*}
for any $n\geq 1$, where
\begin{align} \label{Err-r}
   r_{\g,\theta} (x,\xi) = (2\pi)^{-d} \iint e^{-iy\cdot \eta} \partial_\xi^\g a_1(x,\xi+\theta \eta) D^\g_x a_2(x+y,\xi)dy d\eta,
\end{align}
and $\{r_{\g,\theta}(x,\xi)\}_{|\theta|\leq 1}$ is a bounded symbol of $S^{m_1+m_2 -|\g|}$.
\end{lemma}

Note that, the commutator $i[\Psi_{a}, \Psi_{b}]:= i(\Psi_{a} \Psi_{b} -  \Psi_{b}  \Psi_{a})$
is an operator with symbol in $S^{m_1+m_2-1}$, and the principle symbol is the Poisson bracket
\begin{align*}
    H_{a}b := \{a,b\} = \sum\limits_{j=1}^d \partial_{\xi_j}a \partial_{x_j}b - \partial_{\xi_j}b \partial_{x_j}a.
\end{align*}

The lemma below allows to  estimate  the remainder term \eqref{Err-r}
in the expansion of composition of pseudo-differential operators.

\begin{lemma} \label{Lem-EstiErr}
Let $a, b \in C^\9(\bbr^d\times\bbr^d)$
be such that for any multi-indices $\a, \beta$ with $|\a+\beta|=l$,
$|\p_\xi^\a \p_x^\beta a(x,\xi)| \leq C_1(l) \<x\>^{\rho_1(\beta)} \<\xi\>^{m_1-|\a|}$,
$|\p_\xi^\a \p_x^\beta b(x,\xi)| \leq C_2(l) \<x\>^{\rho_2(\beta)} \<\xi\>^{m_2-|\a|}$,
where $\rho_i(\beta) = \rho_i(|\beta|)$ are decreasing on $|\beta|$.
Set
$$c_\theta(x,\xi)= \iint e^{-i\eta\cdot z} a(x,\xi+\theta\eta) b(x+z,\xi) d\eta dz,$$
where $\theta\in [0,1]$.
Then,
\begin{align} \label{esti-ctheta}
   |c_\theta(x,\xi)|
   \leq C C_1(l) C_2(k) \<x\>^{\rho_1(0)+\rho_2(0)} \<\xi\>^{m_1+m_2}
\end{align}
for $l,k$ such that $l>|\rho_2(0)| +d$,
$k>|m_1|+d $,
where $C$ is independent of $\theta$.
\end{lemma}
(See the Appendix for the proof.)

\begin{corollary} \label{Cor-EstiErr}
Let $a,b$ be as in Lemma \ref{Lem-EstiErr}. For any multi-indices $\a,\a',\beta,\beta'$ and  $\theta\in [0,1]$,
let $c_\theta(x,\xi)$ be as in Lemma \ref{Lem-EstiErr} and
\begin{align*}
  \wt{c}_\theta(x,\xi) := \iint e^{-i\eta\cdot z} \p_\xi^\a \p_x^\beta a(x,\xi+\theta\eta)
                    \p_\xi^{\a'} \p_x^{\beta'} b(x+z,\xi) d\eta dz.
\end{align*}
Then,
there exists $C$ independent of $\theta$ such that
\begin{align*}
     |\p_\xi^\a \p_x^\beta c_\theta(x,\xi)|
     \leq C C_1(|\a+\beta|+l) C_2(|\a+\beta|+k)
          \<x\>^{\rho_1(0)+\rho_2(0)} \<\xi\>^{m_1+m_2-|\a|},
\end{align*}
and
\begin{align*}
    |\wt{c}_\theta(x,\xi)| \leq C C_1(|\a+\beta|+l) C_2(|\a'+\beta'|+k)
    \<x\>^{\rho_1(\beta)+\rho_2(\beta')} \<\xi\>^{m_1+m_2-|\a+\a'|},
\end{align*}
where $C_i(\cdot), \rho_i(\cdot)$, $i=1,2$, $l,k$ are as in Lemma \ref{Lem-EstiErr}.
\end{corollary}

\begin{lemma} \label{Lem-L2-Bdd}
Let  $a\in S^0$, $p\in (1,\9)$.
There exist a constant $C$ and $l \in \mathbb{N}$ such that
\begin{align} \label{pdo-l2}
    \|\Psi_a \|_{\mathcal{L}(L^p)} \leq C |a|_{S^0}^{(l)} .
\end{align}
\end{lemma}

See \cite[Theorem 4.1]{K81} for the case where $p=2$ and
 \cite[Chapter 1.2]{T00} for the general case where $p\in (1,\9)$.

\begin{lemma} \label{Lem-Gar-Ineq}

$(i)$  (G{\aa}rding inequality) Let $a \in S^1$
with ${\rm Re} a(x,\xi) \geq 0$ for $|\xi| \geq R$, $R\geq 1$. Then, there exist
$j_0=j_0(d) \in \mathbb{N}$ and $c(d, R)$, such that
\begin{align} \label{Gard-inequ}
      {\rm Re} \<\Psi_a f, f \> \geq -c(d,R) |a|^{(j_0)}_{S^1} |f|_2^2,\ \ \forall f\in L^2.
\end{align}
Moreover, there exists a constant $C(d)>0$ such that
\begin{align} \label{Gard-R}
  c(d,R) \leq C(d) \<R\>.
\end{align}

$(ii)$ Let $a\in S^{m-1}$, $ m > 2$, ${\rm Re}  a(x,\xi) \geq 0$ for any $|\xi| \geq R$, $R>1$,
and $|\pxi^\a\px^\beta a| \leq C_{\a,\beta} \<x\>^{-2}\<\xi\>^{m-1-|\a|}$ for
any multi-indices $\a, \beta$.
Then,
\begin{align} \label{garding}
   {\rm Re} \<\Psi_a f, f\> \geq - C R \|f\|_{H^{\frac{m-2}{2}}_{-1}}^2,\ \ \forall f\in H^{\frac{m-2}{2}}_{-1},
\end{align}
where $C$ is independent of $R$.
\end{lemma}

{\bf Proof.}
$(i)$ \eqref{Gard-inequ} is the standard Garding estimate, see, e.g., \cite[Lemma 10.3]{LP09}.
As regards \eqref{Gard-R}, let
$\vf$ be a positive smooth function such that $\vf(\xi)=1$ if $|\xi|\leq 1$ and $\vf(\xi)=0$ if $|\xi|\geq 2$.
Set $\vf_R(\xi):= \wt{C}(R) \vf (\frac{\xi}{R})$, where
$\wt{C}(R):=\sup_{x\in\bbr^d,|\xi|\leq R} |a(x,\xi)| \leq |a|_{S^1}^{(0)}\<R\>$.
Then,
${\rm Re}  (a(x,\xi) + \vf_R(\xi)) \geq 0 $, $\forall x, \xi \in \bbr^d$.
By \eqref{Gard-inequ}, there exist $j_0=j_0(d)\in \mathbb{N}$ and a constant $c(d)$, such that
\begin{align} \label{q-vf}
  {\rm Re}  \<\Psi_{a+\vf_R}f, f\>_2
   \geq& -c(d) (|a+\vf_R|^{(j_0)}_{S^1}) |f|_2^2  \nonumber \\
   \geq& -c(d) (|a|^{(j_0)}_{S^1}+|\vf_R|^{(j_0)}_{S^1}) |f|_2^2,\ \ \forall f\in L^2.
\end{align}
Note that, $\vf_R \in S^0$, and for any $l\geq 1$,
$ |\vf_R|_{S^0}^{(l)}
 = \wt{C}(R)|\vf(\frac{\cdot}{R})|_{S^0}^{(l)}
 \leq |a|_{S^1}^{(l)} |\vf|_{S^0}^{(l)} \<R\>$.
Then,
by Lemma \ref{Lem-L2-Bdd},
for some $l\in \mathbb{N}$,
\begin{align} \label{vfR}
    |{\rm Re}  \<\Psi_{\vf_R} f, f\>_2|
    \leq C |a|_{S^1}^{(l)} |\vf|_{S^0}^{(l)} \<R\> |f|_2^2,\ \ \forall f\in L^2.
\end{align}

Moreover, using the facts that  $\<\xi\>^{-1} \leq \<\xi /R\>^{-1}$
and $\<\xi\>^{-1+l} \leq R^{l-1} \<\xi/R\>^{-1+l}$ for any  $ l\geq 1$,
we have $|\vf(\frac{\cdot}{R})|_{S^1}^{(j_0)} \leq |\vf|_{S^1}^{(j_0)}$.
Hence,
\begin{align} \label{vfR*}
     |\vf_R|_{S^1}^{(j_0)} = \wt{C}(R) |\vf(\frac{\cdot}{R})|_{S^1}^{(j_0)}
     \leq  |a|_{S^1}^{(0)} |\vf|_{S^1}^{(j_0)}  \<R\>.
\end{align}
Plugging \eqref{vfR} and \eqref{vfR*} into \eqref{q-vf}, we have
\begin{align*}
   {\rm Re}\<\Psi_a f, f\>_2
   \geq&  - c(d) (|a|^{(j_0)}_{S^1} + |a|^{(0)}_{S^1} |\vf|^{(j_0)}_{S^1} \<R\>) |f|_2^2
       -   C |a|_{S^1}^{(l)} |\vf|_{S^0}^{(l)} \<R\> |f|_2^2
\end{align*}
for some $l\in \mathbb{N}$, which implies \eqref{Gard-R}.

$(ii)$ Let $g=\<x\>^{-1}\<D\>^{\frac{m-2}{2}}f$.
By Lemma \ref{Lem-Err} and Corollary \ref{Cor-EstiErr},
\begin{align*}
     \<\Psi_a f, f\>
   =    \< \<x\>\<D\>^{-\frac{m-2}{2}} \Psi_a \<D\>^{-\frac{m-2}{2}}\<x\> g, g\>
   =   \<(\Psi_{\frac{\<x\>^2 a(x,\xi)}{\<\xi\>^{m-2}}}+ \Psi_r)g, g\>,
\end{align*}
where $r\in S^0$. Note that $ \<\xi\>^{-(m-2)} \<x\>^2 a(x,\xi) \in S^1$,
and $ {\rm Re}  (\<\xi\>^{-(m-2)} \<x\>^2$ $a(x,\xi)) \geq 0$ for any $|\xi|\geq R$.
Then, using \eqref{Gard-inequ} we obtain
\begin{align*}
    {\rm Re} \< \Psi_{\frac{\<x\>^2 a(x,\xi)}{\<\xi\>^{m-2}}} g, g\>
   \geq - C R |g|_2^2
   = - CR \|f\|^2_{H_{-1}^{\frac{m-2}{2}}}.
\end{align*}
Moreover, since $r\in S^0$, by Lemma \ref{Lem-L2-Bdd},
$| {\rm Re} \<\Psi_r g, g\>|
   \leq C|g|_2^2
   = C \|f\|^2_{H_{-1}^{(m-2)/2}}$.
Therefore, combining the estimates above  we obtain \eqref{garding}.
 \hfill $\square$

\begin{lemma} \label{Lem-Inter}
Fix $m>2$. For any $u\in \mathscr{S}$ and any $\ve\in (0,1)$,
\begin{align} \label{inter}
   \|u\|_{H^{\frac{m-2}{2}}_{-1}}
   \leq C  \ve^{\frac 12}  \|u\|_{H^{\frac{m-1}{2}}_{-1}}
        + C  \ve^{-\frac{m-2}{2}}|u|_2,
\end{align}
where $C$ is independent of $\ve$.
In particular,
\begin{align}  \label{inter2}
    \|u\|_{H^{\frac{m-2}{2}}_{-1}}
   \leq C  \|u\|^{\frac{m-2}{m-1}}_{H^{\frac{m-1}{2}}_{-1}}
         |u|^{\frac{1}{m-1}}_2.
\end{align}
\end{lemma}

{\bf Proof.} Let $\theta$ be a smooth  nondecreasing function such that $\theta(\xi)= 0$ for $|\xi|\leq 1$,
and $\theta(\xi) = 1$ for $|\xi|\geq 2$.
Set $\theta_\ve(\xi) := \theta(\ve \xi)$. Then,
\begin{align*}
   \<x\>^{-1}\<D\>^{\frac{m-2}{2}}
   =&\<x\>^{-1}\<D\>^{\frac{m-2}{2}} \theta_\ve(D)
    + \<x\>^{-1}\<D\>^{\frac{m-2}{2}} (1- \theta_\ve(D)) \\
   =&:K_1 + K_2.
\end{align*}
Note that,
$K_1 = a_\ve(x,D) \<x\>^{-1} \<D\>^{(m-1)/2}$,
where $a_\ve(x,D):=  \<x\>^{-1} \theta_\ve(D)$
$\<D\>^{-1/2} \<x\> \in S^0$,
due to Corollary \ref{Cor-EstiErr}.
Since $\<\xi\>^{-1/2} \leq \ve^{1/2}$ on the support of $\theta_\ve$,
$|a_\ve(x,\xi)|_{S^0}^{(l)} \leq \ve^{1/2} C(l)$, $\forall l\geq 1$.
Hence,
\begin{align} \label{k1}
    |K_1 u|_2 \leq \ve^{\frac 12} C(l) |\<x\>^{-1} \<D\>^{\frac{m-1}{2}}u|_2
       =\ve^{\frac 12} C(l) \|u\|_{H^{\frac{m-1}{2}}_{-1}}.
\end{align}
Moreover, since $\<\xi\> \leq 4 \ve^{-1}$ on the support of $1-\theta_\ve(\xi)$,
$|\<\xi\>^{(m-2)/2} (1-\theta_\ve(\xi))|_{S^0}^{(l)}
   \leq C \ve^{-(m-2)/2}$, we have
\begin{align} \label{k2}
   |K_2 u|_2
   \leq |\<D\>^{\frac{m-2}{2}} (1-\theta_\ve(D))u|_2
   \leq C \ve^{-\frac{m-2}{2}} |u|_2.
\end{align}
Thus, \eqref{k1} and \eqref{k2} yield \eqref{inter}.
\eqref{inter2} follows by optimization in $\ve$.  \hfill $\square$

\begin{lemma} \label{Lem-Psi-Esti}
Let $a\in S^0$ and $\theta $ be a smooth nondecreasing function
such that $\theta(\xi) =0$ if $|\xi|\leq 1$ and $\theta(\xi)=1$ if $|\xi|\geq 2$.
Set $\theta_R(\xi) := \theta(\frac \xi R)$,
$c_R (x,\xi):= e^{ M a(x,\xi) \theta_R(\xi)} \in S^0$,
where  $ R, M\geq 1$.
Let $\|a\|_{\9}:= |a|_{L^\9(\bbr^{2d})}$.

$(i)$. There exist  $l\in \mathbb{N}$ and $C(l)>0$ such that for any $R\geq 1$,
\begin{align} \label{esti-Psi-1-1}
       \|\Psi_{c_R}\|_{\mathcal{L}(L^2)} + \|\Psi_{c_R^{-1}}\|_{\mathcal{L}(L^2)}
       \leq C(l) M^l e^{M \|a\|_{\9}}.
\end{align}

$(ii)$ For
$R = C M^{l}e^{2M  \|a\|_{\9}}$
with $C$, $l$ large enough,
$\Psi_{c_R}$ is invertible, and
\begin{align} \label{esti-Psi-inver}
     \|\Psi^{-1}_{c_R}\|_{\mathcal{L}(L^2)}
     \leq C(l)   M^{l} e^{M   \|a\|_{\9}},
\end{align}

$(iii)$
For $ R = C M^{l}e^{2M \|a\|_{\9}}$  with $C$ and $l$ large enough, we have
\begin{align} \label{esti-Psi-h12}
    \|\Psi_{c_R}  \|_{\mathcal{L}(H^{\frac {m-1}{2}}_{-1})}
    + \|\Psi_{c_R}^{-1}  \|_{\mathcal{L}(H^{\frac {m-1}{2}}_{-1})}
    \leq C(l) M^l e^{2M\|a\|_\9}  .
\end{align}

\end{lemma}

The proof of Lemma \ref{Lem-Psi-Esti} is  postponed  to the Appendix.

\section{Homogeneous local smoothing estimates} \label{Sec-Local-Smooth}

This section is mainly concerned with the  local smoothing estimates
for homogeneous solutions to \eqref{equa-rnls}.

\begin{theorem} \label{Thm-Loc-Smooth}
Consider \eqref{equa-rnls} in the homogenous case $f \equiv 0$, i.e.,
\begin{align} \label{equa-hom}
   \partial_t u = i P_t(x,D) u, \ \ u(0) = u_0.
\end{align}
Assume $(A0)$--$(A2)$.
Then, for any $\mathcal{F}_0$-measurable  $u_0\in L^2$, $\bbp$-a.s.,
the solution $u$ to \eqref{equa-rnls}
satisfies $\bbp$-a.s. that
\begin{align} \label{local-smooth-esti-noncons}
     \|u\|_{C([0,T]; L^2)} + \|u\|_{L^2(0,T; H^{\frac{m-1}{2}}_{-1})}
     \leq \wt{D}(e,T) |u_0|_2
\end{align}
Here,
\begin{align} \label{wtD}
  \wt{D}(e,T) = C_1 (1+e_T+|e|_{L^1})^{\frac 12} (1+\beta^*_T+T)^{\frac 12 ((m-2)^2+m)l} e^{C_2(1+\beta^*_T+T)},
\end{align}
for some $l\geq 1$, where $C_1>0$, $C_2\geq 0$,
$\beta^*_T = \sup_{t\in[0,T]} |\beta(t)|$,
and $e$ is as in \eqref{U-l2} below.
In particular, one can take $l=e=1$ and $C_2=0$
in the conservative case.
\end{theorem}

\begin{remark} \label{Rem-Loc}
Similarly to Remark \ref{Rem-Integ}, we have
\begin{align}
    \wt{D}(e,T) \in \bigcap\limits_{1\leq \rho<\9} L^\rho(\Omega),
\end{align}
and when $u_0\in L^\rho(\Omega; L^2)$ for all $\rho \geq 1$,
\begin{align*}
   u\in \bigcap\limits_{1\leq \rho <\9} L^\rho(\Omega; C([0,T];L^2) \cap L^2(0,T; H^{\frac{m-1}{2}}_{-1})).
\end{align*}
\end{remark}

\begin{remark} \label{Rem-Exp-integ}
In the conservative case, for the Schr\"odinger operator ($m=2$) or the
Airy operator ($m=3$),
the constant $\wt{D}(e,T)$ is even exponentially $\bbp$-integrable.
Moreover, if in addition $u_0\in L^\9(\Omega; L^2)$,
then the solution is also exponentially $\bbp$-integrable,
that is, there exists $\delta>0$ such that
\begin{align*}
    \bbe \exp \(\delta (\|u\|_{C([0,T];L^2)} +  \|u\|_{L^2(0,T; H^{\frac{m-1}{2}}_{-1})}) \) < \9.
\end{align*}
\end{remark}

The key role in the proof of Theorem \ref{Thm-Loc-Smooth} is played by the
pseudo-differential operator of order zero constructed in \cite{C02}.
See also \cite{D94,D96,K13,KPV04} in the Schr\"odinger case.

\begin{lemma} (\cite[Lemma 7.1]{C02}) \label{Lem-Doi}
Assume $(A0)$--$(A2)$.
There exist  $\wt{h} (x,\xi)\in S^0$
and $c_1,c_2>0$, such that
\begin{align*}
   H_{\wt{h}} P \leq - c_1\frac{|\xi|^{m-1}}{\<x\>^2} + c_2,
\end{align*}
and $|\p_\xi^\a \p_x^\beta {\wt{h}} (x,\xi)| \leq C_{\a\beta} \<x\>^{-\rho(\beta)} \<\xi\>^{-|\a|}$
for $(x,\xi) \in T^*\bbr^d$, where
$\rho(\beta)$ is equal to $0$ for $\beta=0$,
$1$ for $|\beta|=1$,
$2$ for $|\beta|\geq 2$.
\end{lemma}

As mentioned in Section \ref{Sec-Main},
the Gronwall inequality as in \cite{C02}
will produce a double exponential bound involving Brownian motions,
which,
however, is not $\bbp$-integrable in the non-conservative case.
Instead, we use the martingale property of homogeneous solutions
to control the $L^2$-norm of solutions.
Similar semi-martingales in the dual case
will also be used in the next section.

As in \cite{BRZ14,BRZ16},
we use the notations $U(t,s)$, $s,t\in [0,\9)$,
for evolution operators
corresponding to \eqref{equa-hom}.
Their dual operators  are denoted by $U^*(t,s)$.

\begin{lemma} \label{Lem-U-L2}
$(i)$
For any $\mathcal{F}_0$-measurable  $u_0\in L^2$, $\bbp$-a.s., we have
\begin{align} \label{U-l2}
   |U(t,0) u_0|_2^2
   \leq e(t) |u_0|_2^2,
\end{align}
where $e(t) =   |e^{-\Phi(W(t))}|^2_{L^\9} M(t)$,
\begin{align} \label{Mar}
  M(t) =  \exp\left\{\dd\sum^N_{j=1}\left[\dd\int^t_0 v_j(s)d\b_j(s)-
 \dd\frac12\int^t_0 v^2_j(s)ds\right]\right\},
\end{align}
and $v_j=2{\rm Re}\<\wt{X},\mu_j e_j \wt{X}\>_{2}|\wt{X}|^{-2}_{2}$
with $\wt{X}(t)= e^{\Phi(W(t))} U(t,0)u_0$,  $1\leq j\leq N$.

$(ii)$
For any $\mathcal{F}_0$-measurable $u_0\in L^2$, $\bbp$-a.s., we have
\begin{align} \label{U*-l2}
   |U^*(0,t) u_0|_2^2
   \leq e^*(t) |u_0|_2^2.
\end{align}
Here $e^*(t) = |e^{\Phi(W(t))}|^2_{L^\9} M^*(t)$,
\begin{align} \label{Mar*}
  M^*(t) :=  \exp\left\{\dd\sum^N_{j=1}\left[ - \dd\int^t_0 v^*_j(s)d\b_j(s)-
 \dd \int^t_0 \frac12 (v^*_j)^2(s) ds + \int_0^t \wt{v}_j^*(s)ds\right]\right\},
\end{align}
$v^*_j=2{\rm Re}\<X^*, \mu_j e_j X^*\>_{2} |X^*|^{-2}_{2}$,  and
$\wt{v}_j^* = 2|X^*|^{-2}_{2} {\rm Re} \<X^*, (|\mu_j|^2+\mu^2_j) e_j X^*\>_2$
with $X^*(t) = \ol{e^{-\Phi(W(t))} } U^*(0,t) u_0$, $1\leq j\leq N$.
\end{lemma}

{\bf Proof}
$(i)$. Note that $\wt{X}$ is the homogeneous solution to \eqref{equa-dis}.
Using similar arguments as in the proof of \cite[(1.4)]{BRZ14},
we have that $|\wt{X}|_2^2$ is a continuous martingale with the representation
\begin{align} \label{mar-Y}
    |\wt{X}(t)|_2^2 = |u_0|_2^2 + 2 \sum\limits_{j=1}^N\int_0^t Re \mu_j \<\wt{X}(s),\wt{X}(s)e_j\>_2 d\beta_j(s),\ \  t\in [0,T].
\end{align}
This yields that $|\wt{X}(t)|^2_{2}=|u_0|^2_2 M(t)$, $t\in [0,T]$,
which implies \eqref{U-l2}.

$(ii)$. Since $Id=U(0,t)U(t,0)$,
$ \partial_t U(0,t)=- i U(0,t)\Psi_{P_t}$, we have
\begin{align*}
   \partial_t \<U^*(0,t)u_0,z\>_2
    =&\<u_0,- i U(0,t)\Psi_{P_t} z\>_2
   =\< i \Psi_{P_t}^*U^*(0,t)u_0,z\>, \ \forall z\in H^m,
\end{align*}
which implies that
$v^*(t):= U^*(0,t)u_0$ satisfies the equation
\begin{equation} \label{equ-u*}
  \partial_t v^* = i \Psi_{P_t}^* v^*,\ \ v^*(0)=u_0.
\end{equation}
Then, by It\^o's formula,
\begin{align*}
   d X^* = i P(x,D) X^* dt + \sum\limits_{j=1}^N (\ol{\mu_j}^2 + \frac 12|\mu_j|^2) e_j^2   X^* dt - X^* dW, \ \
   X^*(0)=u_0.
\end{align*}
This yields that
\begin{align} \label{mar-Y*}
    |X^*(t)|_2^2 =& |u_0|_2^2
     + 4  \sum\limits_{j=1}^N \int ({\rm Re} \mu_j )^2  e_j^2  |X^*(t)|^2  dx dt  \nonumber \\
     &- 2 \sum\limits_{j=1}^N\int_0^t {\rm Re} \mu_j \<X^*(s),X^*(s)e_j\>_2 d\beta_j(s),\ \  t\in [0,T].
\end{align}
Thus, $|X^*(t)|^2_{2}=|u|^2_2 M^*(t)$, $t\in [0,T]$,
which implies \eqref{U*-l2}.   \hfill $\square$

Below we first treat the easier conservative case in the spirit of \cite{CKS95}.

{\it \bf Proof of  Theorem \ref{Thm-Loc-Smooth} (Conservative case)}.
Let $\wt{h} \in S^0$ and $\theta$ be as in Lemmas \ref{Lem-Doi} and \ref{Lem-Psi-Esti} respectively.
Set $h(x,\xi):= \wt{h} (x,\xi) \theta(\xi) \in S^0$.
Note that, $\Psi_{P_t}^*=\Psi_{P_t} = e^{-W(t)} \Psi_P e^{W(t)}$ and
\begin{align} \label{cons-l2}
   |u(t)|_2 = |u_0|_2,\ \ t\in [0,T].
\end{align}

Using \eqref{equa-rnls}  we have
\begin{align}\label{partial-u-psip}
  \partial_t {\rm Re}\<u, \Psi_{h} u\>_2
  =& {\rm Re}\<i \Psi_{P_t}u, \Psi_{h} u\>_2 + {\rm Re}\<u, i \Psi_{h} \Psi_{P_t}u \>_2 \nonumber \\
  =& {\rm Re}\<u,i [\Psi_{h}, \Psi_{P_t}]u\>_2.
\end{align}
Since by \eqref{Pt-P},
\begin{align} \label{cons-comm}
    i[\Psi_{h},\Psi_{P_t}]
    = i[\Psi_{h}, \Psi_{P}] + i [\Psi_{h}, e^{-W(t)} [\Psi_{P}, e^{W(t)}]]
    =: \Psi_a + \Psi_{b_t},
\end{align}
where $a(x,\xi) \in S^{m-1}$ and $b_t(x,\xi) \in S^{m-2}$.
It follows that
\begin{align} \label{cons-comm*}
  \partial_t {\rm Re}\<u, \Psi_{h} u\>_2
  = {\rm Re}\<u, \Psi_a u\>_2 + {\rm Re}\<u,\Psi_{b_t}u\>_2.
\end{align}

Below we perform separate estimates of the symbols $a$ and $b_t$
by expanding them up
to zero order via Lemma \ref{Lem-Err}.
The key point here is that the commutator $i[\Phi_h, \Phi_P]$
is an elliptic operator of order $m-1$,
which raises the local regularity of homogeneous solutions,
while the lower order perturbations will then be controlled by the interpolation estimate \eqref{inter}.
We mainly consider the case $m>2$,
the case $m=2$ can be proved similarly.

First, by Lemma \ref{Lem-Err},
\begin{align} \label{cons-a}
      a(x,\xi) =& H_{h} P
          + \sum\limits_{2\leq |\a|\leq m-1} \frac{i}{\a !} (\partial_\xi^\a {h} D^\a_x P - \partial_\xi^\a P D^\a_x {h})
          + r_0(x,\xi) \nonumber \\
        :=& H_{h} P + a_0(x,\xi) + r_0(x,\xi),
\end{align}
where $a_0\in S^{m-2}$ and  $r_0 \in S^0$.
Since by Lemma \ref{Lem-Doi},
$H_{h} P \leq - \frac{c_1}{2} \<x\>^{-2} \<\xi\>^{m-1} + c_2$
for $|\xi|\geq 2$,
using Lemma \ref{Lem-Gar-Ineq} $(ii)$ we get
\begin{align*}
   Re \<(-H_{h} P - \frac{c_1}{2}  \Psi_{\frac{\<\xi\>^{m-1} }{\<x\>^{2}}} + c_2 )u, u \>
   \geq - C \|u\|^2_{H^{\frac{m-2}{2}}_{-1}},
\end{align*}
which implies that
\begin{align} \label{esti-HqP}
   Re \<H_{h} Pu ,u \>
   \leq& -\frac{c_1}{2}  Re \<\Psi_{\frac{\<\xi\>^{m-1}}{\<x\>^{2}}} u , u\> + C \|u\|^2_{H^{\frac{m-2}{2}}_{-1}} + c_2|u|_2^2 \nonumber  \\
   \leq& -\frac{c_1}{2}  \|u\|^2_{H^{\frac{m-1}{2}}_{-1}} + C \|u\|^2_{H^{\frac{m-2}{2}}_{-1}} + C |u|_2^2.
\end{align}
Moreover, by Assumption (A1) and Lemma \ref{Lem-Doi},
$|\p_\xi^\g \p_x^\beta a_0(x,\xi)|\leq C \<x\>^{-2}$ $\<\xi\>^{m-2-|\g|}$,
implying that $\Psi_{\wt{a}_0}:= \<x\>\<D\>^{-\frac{m-2}{2}} \Psi_{a_0}\<D\>^{-\frac{m-2}{2}} \<x\> \in S^0$.
Hence,
\begin{align} \label{esti-b0}
   |\<u, \Psi_{a_0}u\>|
   =&  \bigg|\<\Psi_{\frac{\<\xi\>^{\frac{m-2}{2}}}{\<x\>}}u,
     \Psi_{\wt{a}_0} \Psi_{\frac{\<\xi\>^{\frac{m-2}{2}}}{\<x\>}} u\> \bigg|
   \leq   C \|u\|^2_{H^{\frac{m-2}{2}}_{-1}}.
\end{align}
Thus, taking together \eqref{cons-a}, \eqref{esti-HqP} and \eqref{esti-b0} we obtain
\begin{align} \label{esti-a-cons}
  {\rm Re} \<u,\Psi_au\>
   \leq -\frac{c_1}{2}  \|u\|^2_{H^{\frac{m-1}{2}}_{-1}}
        + C \|u\|^2_{H^{\frac{m-2}{2}}_{-1}}
        + C |u|_2^2.
\end{align}

As regards the symbol $b_t(x,\xi)$,
by Lemma \ref{Lem-Err},
\begin{align} \label{paw*}
      e^{-W(t)}[\Psi_P, e^{W(t)}]
    = \sum\limits_{1\leq |\a|\leq m} \frac{1}{\a !} \Psi_{\partial_\xi^\a P \psi_{\a}(W)},
\end{align}
where $\psi_\a(W) = e^{-W} D_x^\a e^W$,
satisfying that  $|\psi_\a(W)| \leq C |\beta(t)|^{|\a|}$.
Then, applying Lemma \ref{Lem-Err} again we get
\begin{align}  \label{cons-b1r1}
     \Psi_{b_t}
  =& i [\Psi_{h}, \sum\limits_{1\leq|\a|\leq m} \frac{1}{\a !} \Psi_{\partial_\xi^\a P \psi_{\a}(W)}] \nonumber \\
  =& i \sum\limits_{1\leq |\a|\leq m-2} \frac{1}{\a !}
       \bigg[\sum\limits_{1\leq |\beta|\leq m-|\a|-1}
       \frac{1}{\beta !} \big(\partial_\xi^\beta {h} D^\beta_x (\partial^\a_\xi P \psi_\a(W)) \nonumber \\
     & \qquad  - D_x^\beta {h} \partial_\xi^\beta (\partial^\a_\xi P \psi_\a(W)) \big) + \Psi_{r_{2,\a}}\bigg]
       + i \sum\limits_{|\a|=m-1,m}\frac{1}{\a!} [\Psi_{h}, \Psi_{\p_\xi^\a P \psi_\a(W)}] \nonumber \\
    =&: \Psi_{b_{t,1}} + \Psi_{r_{t,1}},
\end{align}
where $\Psi_{r_{t,1}} = i \sum_{1\leq |\a|\leq m-2} \frac{1}{\a !}  \Psi_{r_{2,\a}}
+ i \sum_{|\a|=m-1,m} \frac{1}{\a!} [\Psi_{h}, \Psi_{\p_\xi^\a P \psi_\a(W)}] \in S^0$,
$b_{t,1} = b_t -r_{t,1}$.
Note that, by Assumption $(A1)$,
for any $l\in \mathbb{N}$ and
any multi-indices $\g_1, \g_2$,
\begin{align} \label{esti-r2}
     |r_{t,1}|_{S^0}^{(l)} \leq C(l)(1+ |\beta_t|^{m}),  \ \
     |\p_\xi^{\g_1} \p_x^{\g_2} b_{t,1}(x,\xi)|
     \leq  C (1+|\beta_t|^{m-2}) \frac{\<\xi\>^{m-2-|\g_1|}}{\<x\>^{2}}.
\end{align}
Then, similarly to \eqref{esti-b0}, we have
\begin{align} \label{esti-bt}
  |\<u, \Psi_{b_t} u\>|
  \leq& C (1+|\beta_t|^{m-2}) \|u\|^2_{H^{\frac{m-2}{2}}_{-1}} + C (1+|\beta_t|^m)|u|_2^2.
\end{align}

Now, applying \eqref{esti-a-cons} and \eqref{esti-bt} into \eqref{cons-comm*}
and using \eqref{inter} we get
\begin{align*}
   \p_t {\rm Re} \<u, \Psi_{h} u\>
   \leq& - \frac{c_1}{2}   \|u\|^2_{H^{\frac{m-1}{2}}_{-1}}
        + C(1+|\beta_t|^{m-2}) \|u\|^2_{H^{\frac{m-2}{2}}_{-1}}
        + C(1+|\beta_t|^{m}) |u|^2_2 \\
   \leq& ( - \frac{c_1}{2}  + \ve C (1+|\beta_t|^{m-2}) )  \|u\|^2_{H^{\frac{m-1}{2}}_{-1}}
        + \ve^{-(m-2)} C(1+|\beta_t|^{m}) |u|^2_2.
\end{align*}
Taking  $\ve= c_1(4C(1+|\beta_t|^{m-2}))^{-1}$ we obtain
\begin{align} \label{p-h}
   \partial_t {\rm Re} \<u,\Psi_{h} u\>
   \leq -\frac{c_1}{4} \|u\|^2_{H^{\frac{m-1}{2}}_{-1}}
        + C (1+|\beta_t|^{(m-2)^2+m}) |u|_2^2.
\end{align}
Thus, integrating over $[0,T]$ and using \eqref{cons-l2}  we get
\begin{align*}
      & {\rm Re} \<u(T), \Psi_{h} u(T)\> \nonumber \\
  \leq& {\rm Re} \<u_0, \Psi_{h} u_0\> - \frac{c_1}{4}  \|u\|^2_{L^2(0,T; H^{\frac{m-1}{2}}_{-1})}
        + C T (1+ (\beta^*_T)^{(m-2)^2 +m} ) |u_0|_2^2,
\end{align*}
which implies immediately that
\begin{align}\label{esti-integ-cons}
    \|u\|^2_{L^2(0,T; H^{\frac{m-1}{2}}_{-1})}
    \leq C(1+ T(1+ (\beta^*_T)^{(m-2)^2+m})) |u_0|^2_2,
\end{align}
thereby proving \eqref{local-smooth-esti-noncons} in the  case $m>2$.

The case $m=2$ is  easier. In this case, we do not have the lower order terms $a_0$
and $b_{t,1}$ in \eqref{cons-a} and \eqref{cons-b1r1} respectively.
Hence,  instead  of \eqref{esti-a-cons} and \eqref{esti-bt} we have
\begin{align*}
   {\rm Re} \<u,\Psi_a u\> \leq - \frac{c_1}{2}  \|u\|^2_{H^\frac 12_{-1}} + C|u|_2^2, \ \
   |\<u,\Psi_{b_t} u\>| \leq C(1+|\beta_t|^2) |u|_2^2.
\end{align*}
Therefore, arguing as those below \eqref{esti-bt} we obtain \eqref{esti-integ-cons}
with $m=2$.  \hfill $\square$

Next we treat the non-conservative case,
for which we will use the transformation as in \cite{KPV04}  and
perform the energy method to a new equation.

{\bf Proof of Theorem \ref{Thm-Loc-Smooth} (Non-conservative case)}
Let $\wt{h},\theta$ be as in previous the proof of the conservative case.
Set $\theta_R(\xi):= \theta(\frac{\xi}{R})$,
$h_R(x,\xi) := M \wt{h}(x,\xi) \theta_R(\xi)$, and
$c_R(x,\xi):= \exp\{h_R(x,\xi)\}$ $\in S^0$,
where $M \geq 1$ is to be chosen later, and
\begin{align} \label{R}
   R = C M^l e^{2 M \|\wt{h}\|_\9}
\end{align}
for some $l\geq 1$ and $C$   large enough
such that Lemma \ref{Lem-Psi-Esti} holds.

By \eqref{equa-rnls}, $v:= \Psi_{c_R} u$ satisfies the equation
\begin{align*}
    \pt v = i \Psi_{c_R}\Psi_{P_t} \Psi_{c_R}^{-1} v.
\end{align*}
In view of \eqref{Pt-P} we have
\begin{align*}
  \Psi_{c_R}\Psi_{P_t} \Psi_{c_R}^{-1}
  =& \Psi_P
      + [\Psi_{c_R}, \Psi_P]\Psi^{-1}_{c_R}
       + \Psi_{c_R} e^{-\Phi(W)}[\Psi_{P},  e^{\Phi(W)}] \Psi^{-1}_{c_R}.
\end{align*}
Thus,
\begin{align} \label{noncons-dl2}
      \frac12 \pt|v|_2^2
      =& Re \<v, i [\Psi_{c_R}, \Psi_{P}]\Psi^{-1}_{c_R} v\>
       + Re\<v, i \Psi_{c_R} e^{-\Phi(W)}[ \Psi_{P},  e^{\Phi(W)}] \Psi_{c_R}^{-1} v\> \nonumber \\
      =&: Re \<v, \Psi_a v\> + Re \<v, \Psi_{b_t} v\>,
\end{align}
where $\Psi_a= i [\Psi_{c_R}, \Psi_{P}]\Psi^{-1}_{c_R}$,
$\Psi_{b_t} = i \Psi_{c_R} e^{-\Phi(W)}[ \Psi_{P},  e^{\Phi(W)}] \Psi_{c_R}^{-1}$.

Note that, unlike the conservative case, the symbols
$a, b_t$ have the same order $m-1$.
Below we separate, via Lemma \ref{Lem-Err},
the principle symbols and the lower order symbols of $a$ and $b_t$.
As in the conservative case,
we will mainly consider the case $m>2$,
the case $m=2$ can be proved similarly.

First, for the symbol $a$, note that
\begin{align*}
  i [ \Psi_{c_R},  \Psi_P]
  =   \Psi_{H_{c_R} P} + \Psi_{r_0'}
  =  \Psi_{H_{h_R} P} \Psi_{c_R} + \Psi_{r_0''} + \Psi_{r_0'},
\end{align*}
where $r_0', r_0'' \in S^{m-2}$.
This implies that
\begin{align} \label{sym-a}
   \Psi_a
   = \Psi_{H_{h_R} P} + \Psi_{r_0},
\end{align}
where $\Psi_{r_0} = (\Psi_{r_0'}+ \Psi_{r_0''}) \Psi_{c_R}^{-1}$.
By Assumption $(A1)$, Corollary \ref{Cor-EstiErr} and straightforward computations, we have that
for
any multi-indices $\a,\beta$, $|\a+\beta|=l$,
there exists $l'\geq 1$ such that
\begin{align} \label{esti-r1'}
  |\p_\xi^\a \p_x^\beta r_0 |
  \leq C(l') M^{l'} e^{2M\|\wt{h}\|_{\9}} \frac{\<\xi\>^{m-2 - |\a|}}{\<x\>^2}.
\end{align}

As regards the symbol $b_t(x,\xi)$, we compute
\begin{align} \label{deco-bt}
   \Psi_{b_t}
   =& i e^{-\Phi(W)} [\Psi_P, e^{\Phi(W)}] + i [\Psi_{c_R}, e^{-\Phi(W)}[\Psi_P, e^{\Phi(W)}]] \Psi^{-1}_{c_R} \nonumber \\
   =& \sum\limits_{|\a|=1} i  \partial^\a_\xi P \psi_\a(\Phi(W)) \nonumber \\
   &   + \bigg(\sum\limits_{2\leq |\a|\leq m} \frac{i}{\a !}  \partial^\a_\xi P \psi_\a(\Phi(W))
       + i [\Psi_{c_R}, e^{-\Phi(W)}[\Psi_P, e^{\Phi(W)}]] \Psi^{-1}_{c_R} \bigg) \nonumber \\
   =:& \Psi_{b_{t,1}} + \Psi_{b_{t,2}},
\end{align}
where $\psi_\a(\Phi(W)) = e^{-\Phi(W)} D^\a_x e^{\Phi(W)}$.
Note that, by Assumption $(A1)$,
\begin{align} \label{bt1}
    |b_{t,1}(x,\xi)| \leq C(1+|\beta_t|+t) \<x\>^{-2} \<\xi\>^{m-1},
\end{align}
and for any multi-indices $\g, \delta$, $|\g+\delta|=l$,
there exists some $l'$ such that
\begin{align} \label{esti-b2}
  |\p_\xi^{\g} \p_x^{\delta} b_{t,2}|
  \leq C(l') (1+|\beta(t)|+t)^{m} M^{ l'} e^{2M\|\wt{h}\|_{\9}} \frac{\<\xi\>^{m-2-|\g|}}{\<x\>^2}.
\end{align}

Thus, it follows from \eqref{noncons-dl2}, \eqref{sym-a} and \eqref{deco-bt} that
\begin{align} \label{esti-b2*}
   \frac 12 \p_t |v|_2^2
   = {\rm Re} \<v, (\Psi_{H_{h_R}P}+ \Psi_{b_{t,1}})v\>
   + {\rm Re} \<v,(\Psi_{r_0} + \Psi_{b_{t,2}}) v\>.
\end{align}

For the first term in the right hand side of \eqref{esti-b2*},
taking into account Lemma \ref{Lem-Doi}
and the boundedness of $b_{t,1}$, we take
\begin{align} \label{M}
 M= \frac{4C}{c_1} (1+\beta^*_T+T)
\end{align}
and get
\begin{align} \label{ReM}
    \frac 1 M {\rm Re} (H_{h_R}P+b_{t,1})(x,\xi) \leq - \frac {c_1}{4} \frac{\<\xi\>^{m-1}}{\<x\>^2} + \frac{c_2}{M}, \ \ |\xi|\geq 2R.
\end{align}
Then, arguing as in the proof of \eqref{esti-HqP} we obtain
\begin{align} \label{M*}
    {\rm Re} \<v, (\Psi_{H_{h_R}P} + \Psi_{b_{t,1}})v\>
   \leq & - \frac {c_1}{4} M \|v\|^2_{H^{\frac{m-1}{2}}_{-1}} + C R M \|v\|^2_{H^{\frac{m-2}{2}}_{-1}} + C|v|_2^2.
\end{align}

Regarding the second term in the right hand side of \eqref{esti-b2*},
similarly to \eqref{esti-b0},
using \eqref{esti-r1'}, \eqref{esti-b2}
we get
\begin{align} \label{M**}
  Re \<v, (\Psi_{r_0} +\Psi_{b_{t,2}})v\>
  \leq& C(l') (1+|\beta_t|+t)^m M^{l'} e^{2M\|\wt{h}\|_{\9}} \|v\|^2_{H^{\frac{m-2}{2}}_{-1}}.
\end{align}

Now, plugging \eqref{M*} and \eqref{M**} into  \eqref{esti-b2*}
and using \eqref{inter} we get
\begin{align*}
  \frac 12 \p_t|v|_2^2
  \leq& (-\frac {c_1}{4} M + C(1+|\beta_t|+t)^m M^{l''}e^{2M\|\wt{h}\|_{\9}}\ve) \|v\|^2_{H^{\frac{m-1}{2}}_{-1}} \\
      & + C(1+|\beta_t|+t)^m M^{l''} e^{2M\|\wt{h}\|_{\9}} \ve^{-(m-2)} |v|_2^2
\end{align*}
for some $l''\geq 1$. Choosing $\ve= c_1 (8C)^{-1} (1+|\beta_t|+t)^{-m} M^{1-l''} e^{-2M\|\wt{h}\|_{\9}}$ yields
\begin{align} \label{dv}
   \frac 12 \p_t |v|_2^2
   \leq -\frac {c_1}{8} M \|v\|^2_{H^{\frac{m-1}{2}}_{-1}}
        + C (1+|\beta_t|+t)^{(m-1)m} M^{l''(m-1)} e^{CM} |v|_2^2,
\end{align}
which implies that
\begin{align} \label{esti-nonc}
   \frac {c_1}{8}M  \|v\|^2_{L^2(0,T; H^{\frac{m-2}{2}}_{-1})}
  \leq& \frac 12 (|v_0|_2^2 +|v(T)|_2^2)  \\
      &  + C(1+\beta^*_T+T)^{(m-1)m} M^{l''(m-1)} e^{CM} \|v\|^2_{L^2(0,T; L^2)}.  \nonumber
\end{align}

Note that, applying Gronwall's inequality to \eqref{dv} implies that
\begin{align} \label{d-exp}
   \|v(t)\|_{L^\9(0,T;L^2)}^2 \leq C|v(0)|_2^2 \exp \{C(1+\beta^*_T+T)^{(m-1)m} M^{l''(m-1)} e^{CM}\},
\end{align}
which includes a non-integrable double exponential of Brownian motions
and so can not yield, via \eqref{esti-nonc}, the integrability of the  $L^2(0,T; H^{(m-1)/2}_{-1})$-norm of $v$.

Instead, we use the martingale property of homogeneous solutions in Lemma \ref{Lem-U-L2} (i).
Then, taking into account the
boundedness of $\Psi_{c_R}, \Psi^{-1}_{c_R}$
in   $H^{(m-1)/2}_{-1}$ and $L^2$, we obtain \eqref{local-smooth-esti-noncons}
for $m>2$ in the nonconservative case.

The case $m=2$ is easier. In this case, similarly to \eqref{M*},
applying Lemma \ref{Lem-Gar-Ineq} $(i)$ we have
\begin{align*}
   {\rm Re} \<(H_{h_R}P + b_{t,1})v, v\>
   \leq -\frac{c_1}{4} M \|v\|^2_{H^{\frac 12}_{-1}} + CRM |v|_2^2.
\end{align*}
Moreover, we have that $r_0+ b_{t,2} \in S^0$,
which implies \eqref{M**}
with $H^{(m-2)/2}_{-1}$ replaced by $L^2$.
Then,
similar arguments as above yield \eqref{esti-nonc} with $m=2$,
thereby proving \eqref{local-smooth-esti-noncons}.
The proof is complete. \hfill $\square$

\section{Proof of main results} \label{Sec-Stri-Esti}

We start with the Strichartz  and local smoothing estimates of the
free  group $\{e^{-itP(x,D)}\}$. For simplicity,
we use the abbreviations $\calx_{s,p,q}$, $\calx'_{-s,p',q'}$
for $\calx_{T,s,p,q}$ and $\calx'_{T,-s,p',q'}$ respectively.

\begin{lemma}  \label{Lem-stsm-P}
Assume $(A0)$--$(A3)$. For any $u_0\in L^2 $ and admissible triple $(s,p,q) \in \mathcal{A}$,
\begin{align} \label{homo-stsm-P}
    \|e^{i\cdot P(x,D)} u_0\|_{\calx_{s,p,q}} \leq C (1+T)^{\frac 12} |u_0|_2.
\end{align}
Moreover, for any  $(s_1, p_1, q_1), (s_2, p_2, q_2) \in \mathcal{A}$ and any $f\in \calx'_{-s_2, p_2',q_2'}$,
\begin{align} \label{inho-stsm-P}
  \bigg\|\int_0^\cdot e^{ i(\cdot-s) P(x,D)} f(s)ds \bigg\|_{\calx_{s_1, p_1,q_1}}
  \leq C (1+T) \|f\|_{\calx'_{-s_2, p'_2,q'_2}},
\end{align}
where  $C$ is independent of $T$.
\end{lemma}

{\bf Proof.}
In order to prove \eqref{homo-stsm-P},
in view of Assumption {\rm $(A3)$},
we only need to prove  that
\begin{align} \label{local-P}
   \|e^{i\cdot P(x,D)} u_0\|^2_{L^2(0,T; H^{\frac{m-1}{2}}_{-1})}
   \leq C(1+T) |u_0|_2^2.
\end{align}

For this purpose, let $u(t)= e^{it P(x,D)}u_0$
and $h$ be as in the proof of Theorem \ref{Thm-Loc-Smooth}.
Similarly to \eqref{esti-a-cons},
\begin{align*}
   \p_t {\rm Re} \<u, \Psi_{h} u\>
   \leq - \frac{c_1}{2} \|u\|^2_{H^{\frac{m-1}{2}}_{-1}} + C \|u\|^2_{H^{\frac{m-2}{2}}_{-1}} + C|u|_2^2,
\end{align*}
which along with \eqref{inter} implies that
\begin{align*}
   \p_t {\rm Re} \<u, \Psi_{h} u\>
   \leq (- \frac{c_1}{2} + C\ve) \|u\|^2_{H^{\frac{m-1}{2}}_{-1}} + C(1+\ve^{-(m-2)})|u|_2^2.
\end{align*}
Thus, taking $\ve$ small enough we obtain \eqref{local-P},
thereby proving \eqref{homo-stsm-P}.

Regarding \eqref{inho-stsm-P}, we
note that,
for any $(s,p,q) \in \mathcal{A}$ and $z\in L^2$,
by \eqref{homo-stsm-P},
\begin{align*}
   \<\int_0^T e^{isP(x,D)} f(s)ds, z\>_2
   =&  \int_0^T\< f(s), e^{-isP(x,D)} z\>_2ds \nonumber\\
   \leq&  C(1+T)^\frac 12 \| f\|_{\calx'_{-s, p',q'}} |z|_2,
\end{align*}
which implies that
\begin{align} \label{esti-p-inT}
   \bigg|\int_0^T e^{-i s P(x,D)} f(s) ds \bigg|_2 \leq C (1+T)^\frac 12 \| f\|_{\calx'_{-s, p',q'}}.
\end{align}

Now, let $f= f_1 + f_2$, $f_1 \in L^2(0,T; H^{-(m-1)/2}_{1})$ and
$f_2 \in L^{q_2'}(0,T; \dot{W}^{-s_2,p_2'})$.
Since
$\int_0^T e^{i(t-s)P(x,D)}f_1(s)ds = e^{i t P(x,D)} \int_0^T e^{ -is P(x,D)}f_1(s)ds$,
by \eqref{homo-stsm-P} and \eqref{esti-p-inT},
\begin{align*}
     \bigg\|\int_0^T e^{i(\cdot-s)P(x,D)}f_1(s)ds \bigg\|_{L^{q_1}(0,T;\dot{W}^{s_1,p_1})}
   \leq& C (1+T)^\frac 12 \bigg| \int_0^T e^{ -is P(x,D)}f_1(s)ds \bigg|_2\nonumber \\
   \leq& C(1+T )  \| f_1\|_{L^2(0,T; H^{-\frac{m-1}{2}}_{1})}.
\end{align*}
Then, since $q_1>2$, by Christ-Kiselev's lemma (see e.g. \cite[Lemma 3.1]{SS00}
\footnote{The proof of \cite[Lemma 3.1]{SS00} also works for homogeneous Sobolev spaces.}),
\begin{align} \label{esti-f1}
     \bigg\|\int_0^\cdot e^{i(\cdot-s)P(x,D)}f_1(s)ds \bigg\|_{L^{q_1}(0,T; \dot{W}^{s_1,p_1})}
   \leq C (1+T) \| f_1 \|_{L^2(0,T; H^{-\frac{m-1}{2}}_{1})}.
\end{align}

Similarly, since $q_1>q_2'$, similar arguments as above imply that
\begin{align} \label{esti-f2}
     \bigg\|\int_0^\cdot e^{i(\cdot-s)P(x,D)}f_2(s)ds \bigg\|_{L^{q_1}(0,T; \dot{W}^{s_1,p_1})}
   \leq C(1+T)   \| f_2 \|_{L^{q_2'}(0,T; \dot{W}^{-s_2,p_2'})}.
\end{align}
Thus, setting $v(t):= \int_0^t e^{-i(t-s)P(x,D)}f ds $ we obtain
\begin{align} \label{esti-lqhp}
     \|v\|_{L^{q_1}(0,T; \dot{W}^{s_1,p_1})}
   \leq C (1+T)  \| f \|_{\calx'_{-s_2, p_2',q_2'}}.
\end{align}

As regards the estimate for the $L^2(0,T; H^{(m-1)/2}_{-1})$-norm of $v$,
we first note that $\p_t|v(t)|_2^2 =2 {\rm Re} \int v(t,x) \ol{f(t,x)} dx$,
which implies that
\begin{align} \label{v-l2-con}
    \|v\|^2_{L^\9(0,T; L^2)}
     \leq  \|v \|_{\calx_{s_2,p_2,q_2}}   \|f\|_{\calx'_{-s_2, p_2',q_2'}}.
\end{align}
Moreover, arguing as in the proof of \eqref{esti-integ-cons}
we have
\begin{align*}
   \|v \|^2_{L^2(0,T; H^{\frac{m-1}{2}}_{-1}) }
   \leq& C|v(T)|_2^2  + C\|v\|^2_{L^2(0,T; L^2)} + C    |\<f, \Psi_{h} v\>| + |\<v, \Psi_{h} f\>| ,
\end{align*}
where $\<\ ,\ \>$ is the pairing between $\calx_{s_2,p_2,q_2}$ and $\calx'_{-s_2,p'_2,q'_2}$,
Thus, by \eqref{v-l2-con},
\begin{align*}
   \|v \|^2_{L^2(0,T; H^{\frac{m-1}{2}}_{-1}) }
   \leq&  C (1+T) \|f\|_{\calx'_{-s_2,p_2',q_2'}} \|v\|_{\calx_{s_2,p_2,q_2}}  \nonumber \\
   \leq& \frac 12 \|v \|^2_{\calx_{s_2,p_2,q_2}} + C^2 (1+T)^2 \|f\|^2_{\calx'_{-s_2, p_2',q_2'}},
\end{align*}
which along with \eqref{esti-lqhp} implies \eqref{inho-stsm-P}.
The proof is complete. \hfill $\square$

Below we prove Strichartz and local smoothing estimates for the homogeneous solutions to \eqref{equa-hom}.
\begin{lemma} \label{Homo-U}
Assume $(A0)$--$(A3)$. For each $\mathcal{F}_0$-measurable $u_0\in L^2$, $\bbp$-a.s.,
and $(s,p,q)\in \mathcal{A}$, we have $\bbp$-a.s.
\begin{equation} \label{homo-stri}
    \|U(\cdot,0)u_0\|_{\calx_{s,p,q}} \leq  C(e,T)   |u_0|_2,
\end{equation}
where $C(e,T):=C(1+T)  (1+\beta^*_T +T)^m \wt{D}(e,T)$
with $\wt{D}(e,T)$ as  in \eqref{local-smooth-esti-noncons}.
\end{lemma}

{\bf Proof.} Set $u(t):= U(t,0)u_0$,
$b(t,x,D) : =  i e^{-\Phi(W)(t,x)}  [\Psi_{P},  e^{\Phi(W)(t,x)} ]$, $t\in [0,T]$.
By \eqref{equa-hom},
\begin{align*}
   u(t)    = e^{itP(x,D)} u_0 + \int_0^t e^{i(t-s)P(x,D)} b(s,x,D)u(s) ds.
\end{align*}
Then, by Lemma \ref{Lem-stsm-P},
\begin{align*}
    \| u \|_{\calx_{s,p,q}}
    \leq C (1+T) (|u_0|_2 + \|b(\cdot,\cdot,D)u\|_{L^2(0,T; H^{-\frac{m-1}{2}}_{1})}).
\end{align*}
Note that, by Corollary \ref{Cor-EstiErr}, the $S_0^{(l)}$ semi-norm  of
$\<x\>\<D\>^{-(m-1)/2} b(t,x,D)$ $\<D\>^{-(m-1)/2}\<x\>$
is bounded by $C(l') (1+|\beta_t|+t)^m$ for some $l'\geq 1$.
Then,
\begin{align} \label{bu}
   \|b(\cdot,\cdot,D)u\|_{L^2(0,T;H^{-\frac{m-1}{2}}_{1})}
   \leq C (1+\beta^*_T +T)^m   \|u\|_{L^2(0,T; H^{\frac{m-1}{2}}_{-1})}.
\end{align}
Therefore, in view of Theorem \ref{Thm-Loc-Smooth},
we obtain \eqref{homo-stri}. \hfill $\square$

Similarly, for the dual operator $U^*(0,t)$ we have
\begin{lemma} \label{Homo-U*}
Assume $(A0)$--$(A3)$.  For each $\mathcal{F}_0$-measurable $u_0\in L^2$, $\bbp$-a.s.,
and for any $(s,p,q) \in \mathcal{A}$,
\begin{equation} \label{homo-stri-U*}
    \|U^*(0,\cdot)u_0\|_{\calx_{s,p,q}} \leq  C(e^*,T) |u_0|_2,\ \ \bbp-a.s.,
\end{equation}
where $ C(e^*,T)$ is as in \eqref{homo-stri} with $e^*$ replacing $e$.
\end{lemma}

{\bf Proof.}
We note that,
\begin{align*}
  P^*_t
  =&  \ol{e^{\Phi(W)}}  P  \ol{e^{ - \Phi(W)}}
  =  P + \ol{e^{\Phi(W)}}[P, \ol{e^{ - \Phi(W)}}]   \\
  =&  P + \sum\limits_{1\leq |\a|\leq m} \frac{1}{\a!} \p_\xi^\a P \ol{\psi_\a(-\Phi(W))}
\end{align*}
with $\psi_\a(-\Phi(W))= e^{\Phi(W)} D_x^\a e^{-\Phi(W)}$,
which has similar structure as $P_t$.
Hence, arguing as in the proof of Theorem \ref{Thm-Loc-Smooth}
and using \eqref{U*-l2},
we also have the homogeneous local smoothing estimates for $U^*(0,t)u_0$,
which consequently yields \eqref{homo-stri-U*}
by similar arguments as in the previous proof of Lemma \ref{Homo-U}.  \hfill $\square$

\begin{corollary} \label{Dual-Esti}
Assume $(A0)$--$(A3)$.  For any $(s,p,q)\in \mathcal{A}$ and
$\{\mathcal{F}_t\}$-adapted $f\in \calx'_{-s,p',q'}$, $\bbp$-a.s.,
\begin{equation} \label{esti-f-dual}
   \bigg|\int_0^T U(0,s)f(s)ds \bigg|_2 \leq   C(e^*,T) \|f\|_{\calx'_{-s,p',q'}},\ \ \bbp-a.s.,
\end{equation}
where $ C(e^*,T)$ is as in \eqref{homo-stri-U*}.
\end{corollary}

{\it Proof.}
By Lemma \ref{Homo-U*}, for any $z\in L^2$,
\begin{align*}
   \<\int_0^T U(0,s)f(s)ds, z\>_2
   =& \int_0^T \< f(s), U^*(0,s)z\>_2 ds
   \leq \|f\|_{\calx'_{-s,p',q'}} \|U^*(0,\cdot)z\|_{\calx_{s,p,q}}  \nonumber \\
   \leq&   C(e^*,T)\|f\|_{\calx'_{-s,p',q'}}  |z|_2,
\end{align*}
which implies \eqref{esti-f-dual} by duality.
\hfill $\square$ \\

{\it \bf Proof of Theorem \ref{Thm-Stri-RNLS}.}
$(i)$. We reformulate
(\ref{equa-rnls}) in the mild form
$$u(t)=U(t,0) u_0 + \int_0^t U(t,s)f(s)ds $$
with $f= e^{-\Phi(W)} F$.
By Lemma \ref{Homo-U},
we only need to prove \eqref{l-stri-esti-lps} when $u_0 \equiv 0$.

First we prove that for
any  $(s_i,p_i,q_i) \in \mathcal{A}$, $i=1,2$,
\begin{align} \label{Lpq-Xpq}
   \|u\|_{L^{q_1}(0,T; \dot{W}^{s_1,p_1})}
    \leq   C(e,T) C(e^*,T)  \|f\|_{\calx'_{-s_2,p_2',q_2'}}.
\end{align}
In particular, taking $(s_1,p_1,q_1)=(0,2,\9)$
we have that for any  $(s,p,q)\in \mathcal{A}$,
\begin{align} \label{Lpq-Xpq*}
   \|u\|_{L^\9(0,T; L^2)} \leq C(e,T) C(e^*,T) \|f\|_{\calx'_{-s, p',q'}}.
\end{align}

The arguments is similar to that of \eqref{esti-lqhp}.
In fact, let $f= f_1 + f_2$, $f_1\in L^{q_2'}(0,T; \dot{W}^{-s_2,p_2'})$ and $f_2\in L^2(0,T; H^{-\frac {m-1}{2}}_{1})$.
Note that
\begin{align*}
    \int_0^T U(t,s) f_1(s) ds
   = \int_0^T U(t,0)U(0,s) f_1(s) ds
   = U(t,0) \int_0^T U(0,s) f_1(s) ds,
\end{align*}
which along with Lemma \ref{Homo-U}  and Corollary \ref{Dual-Esti} implies that
\begin{align*}
       \bigg\|\int_0^T U(\cdot,s) f_1(s) ds \bigg\|_{L^{q_1}(0,T; \dot{W}^{s_1,p_1})}
    \leq&   C (e,T ) \bigg|\int_0^T U(0,s) f_1(s) ds \bigg|_2 \nonumber \\
    \leq&   C (e,T)C(e^*,T)  \|f_1\|_{L^{2}(0,T; H^{-\frac{m-1}{2}}_1)}.
\end{align*}
Since $q_1> 2$, using the Christ-Kiselev lemma  we obtain
\begin{align*}
   \bigg \| \int_0^\cdot U(\cdot,s) f_1(s) ds \bigg\|_{L^{q_1}(0,T; \dot{W}^{s_1,p_1})}
    \leq    C (e,T)C(e^*,T)   \|f_1\|_{L^{2}(0,T;  H^{-\frac{m-1}{2}}_1)}.
\end{align*}
Similarly,
since $q_1>q'_2$,
similar arguments as above yield that
\begin{align*}
   \bigg \| \int_0^\cdot U(\cdot,s) f_2(s) ds \bigg\|_{L^{q_1}(0,T; \dot{W}^{s_1,p_1})}
    \leq    C (e,T)C(e^*,T)   \|f_2\|_{L^{q'_2}(0,T; \dot{W}^{-s_2,p'_2})}.
\end{align*}
Thus,  combining these estimates together we obtain \eqref{Lpq-Xpq}, as claimed.

Below we estimate the $L^2(0,T; H^{(m-1)/2}_{-1})$-norm of $u$.
We shall consider the conservative and non-conservative cases respectively.

{\it Conservative case}.
Similarly to \eqref{esti-integ-cons}, we have
\begin{align} \label{esti-H12-F}
  \|u\|^2_{L^2(0,T; H^{\frac {m-1}{2}}_{-1})}
  \leq& C (1+ T(1+(\beta^*_T)^{(m-2)^2+m})) \|u\|_{L^\9(0,T; L^2)}^2 \nonumber \\
      & + C (|  \<f,\Psi_{h} u\>| +  |\<f, \Psi^*_{h} u\> |),
\end{align}
where $\<\ ,\ \>$ is the pairing between $\calx_{s_2,p_2,q_2}$ and $\calx'_{-s_2,p'_2,q'_2}$,
and ${h} \in S^0$ is the symbol as in the proof of Theorem \ref{Thm-Loc-Smooth}.
Since  $\Psi_{h}, \Psi^*_{h} \in \mathcal{L}(\calx_{s_2,p_2,q_2})$,
the last term of the right hand side above is bounded by
\begin{align*}
    C \|f\|_{\calx_{-s_2,p'_2,q'_2}}  \|u\|_{\calx_{s_2,p_2,q_2}}
   \leq& \frac 14 \|u\|^2_{\calx_{s_2,p_2,q_2}} + 4 C  \|f\|^2_{\calx_{-s_2,p'_2,q'_2}}.
\end{align*}
Plugging this into \eqref{esti-H12-F} and using \eqref{Lpq-Xpq*} we get
\begin{align*}
    \|u\|^2_{L^2(0,T; H^{\frac {m-1}{2}}_{-1})}
    \leq& \frac 1 4  \|u\|^2_{\calx_{s_2,p_2,q_2}} \\
        & + C  (1+T(1+ |\beta^*_T|^{(m-2)^2+m})) (C(e,T)C(e^*,T))^2 \|f\|^2_{\calx'_{-s_2,p'_2,q'_2}},
\end{align*}
which along with \eqref{Lpq-Xpq} implies that
\begin{align*}
   \|u\|^2_{\mathcal{X}_{s_1,p_1,q_1}}
   \leq C  (1+ T(1+ (\beta^*_T)^{(m-2)^2+m})) (C(e,T) C(e^*,T))^2 \|f\|^2_{\calx'_{-s_2,p'_2,q'_2}},
\end{align*}
thereby proving   \eqref{l-stri-esti-lps} in the conservative case.

{\it Non-conservative case.}
Let $v: = \Psi_{c_R}u$, where $c_R$ is as in the proof of Theorem \ref{Thm-Loc-Smooth}.
Then, $v$ satisfies the equation
\begin{align*}
  \partial_t v =  i \Psi_{c_R} \Psi_{P_t} \Psi^{-1}_{c_R} v + \Psi_{c_R} f,\ \ v(0)=0.
\end{align*}
Arguing as in the proof of \eqref{esti-nonc} and using \eqref{Lpq-Xpq*} we get
\begin{align} \label{v-f}
         \|v\|^2_{L^2(0,T; H^{\frac {m-1}{2}}_{-1})}
   \leq& C  |v_T|_2^2
         + C(1+\beta^*_T+T)^{(m-1)m} M^{l(m-1)} e^{CM} \|v\|^2_{L^2(0,T; L^2)} \nonumber \\
       &  + C |\<v, \Psi_{c_R}f\>|
\end{align}
for some $l \geq 1$, where $M= 4c_1^{-1}C (1+ \beta^*_T+T)$.
Since $\Psi^{-1}_{c_R} \in \mathcal{L}(H^{\frac{m-1}{2}}_{-1})$,
$\Psi_{c_R} \in \mathcal{L}(L^2) \cap \mathcal{L}(\calx_{s,p,q})$,
$\Psi^*_{c_R} \in \mathcal{L}(\calx_{s,p,q})$,
and the  norms are bounded by the semi-norms of
$c_R^{-1}$ and $c_R$, we get
\begin{align*}
      & \|u\|^2_{L^2(0,T;H^{\frac {m-1}{2}}_{-1})}  \nonumber \\
  \leq& C M^{l'} e^{CM} (1+T) (1+\beta^*_T+T)^{(m-1)m} (C(e,T)C(e^*,T))^2 \|f\|^2_{\calx'_{-s_2,p'_2,q'_2}}   \nonumber \\
      & + C M^{l'} e^{CM} \|u\|_{\calx_{s_2,p_2,q_2}} \|f\|_{\calx'_{-s_2,p'_2,q'_2}} \nonumber \\
  \leq&    C M^{l'} e^{CM} (1+T) (1+\beta^*_T+T)^{(m-1)m} (C(e,T) C(e^*,T))^2 \|f\|^2_{\calx'_{-s_2,p'_2,q'_2}}
       + \frac 12 \|u\|^2_{\calx_{s_2,p_2,q_2}},
\end{align*}
which along with \eqref{Lpq-Xpq} implies \eqref{l-stri-esti-lps}.
The statement $(i)$ is  proved.

$(ii)$. For each $1\leq j\leq d$, let $w_j :=\p_{x_j}u$.
By \eqref{equa-rnls},
\begin{align*}
   \p_t w_j
   = i P_t(x,D) w_j + i [\p_{x_j}, P_t(x,D)]u + \p_{x_j} f.
\end{align*}
Let $g(t,x,D)= i[\p_{x_j}, P_t(x,D)]$.
Similarly to \eqref{paw*},
\begin{align*}
   P_t = P + e^{-\Phi(W)} [P,e^{\Phi(W)}]
   = P + \sum\limits_{1\leq |\a|\leq m}
     \frac{1}{\a !} \Psi_{\p_\xi^\a P \psi_\a(\Phi(W))}.
\end{align*}
Then, since $P$ is independent of $x$, we have
$$
 g(t,x,\xi)
   = (-1)  \sum_{1\leq |\a|\leq m}  \frac{1}{\a !} \p_\xi^\a P  D_{x_j}\psi_\a(\Phi(W)),
$$
which implies that
\begin{align*}
   |\p_{\xi}^\beta \p_x^\g g(t,x,\xi)|
   \leq C_{\beta \g} (1+\beta^*_T + T)^m \frac{\<\xi\>^{m-1-|\beta|}}{\<x\>^2}
\end{align*}
for any multi-indices $\beta$, $\g$. Thus, applying \eqref{Stri-X} we obtain
\begin{align*}
       & \|w_j\|_{\mathcal{X}_{T,s_1,p_1,q_1}} \\
   \leq& D(e,e^*,T) (|\p_{x_j} u(0)|_2 + \|g(t,x,D) u\|_{L^2(0,T; H^{-\frac{m-1}{2}}_1)}
        + \|\p_{x_j}f \|_{\mathcal{X}'_{T,-s_3,p_3',q_3'}}   ) \\
   \leq& C D(e,e^*,T) (|u(0)|_{H^1} +(1+\beta^*_T + T)^m  \| u\|_{L^2(0,T; H^{\frac{m-1}{2}}_{-1})}
         + \|\p_{x_j}f \|_{\mathcal{X}'_{T,-s_3,p_3',q_3'}}   ),
\end{align*}
which implies \eqref{l-stri-esti-wps}. The proof of Theorem \ref{Thm-Stri-RNLS} is complete. \hfill $\square$\\

{\bf Proof of Theorem \ref{Thm-Stri-per}.}
The proof follows the lines as those in the proof of Theorems \ref{Thm-Loc-Smooth} and \ref{Thm-Stri-RNLS}.
In fact,
we can   derive similar equations as \eqref{partial-u-psip} and \eqref{noncons-dl2}
in the conservative and non-conservative cases respectively.
In each case,
the lower perturbation only contributes the $H^{(m-2)/2}_{-1}$-norm
of $u$,
which can be controlled by the interpolation estimate \eqref{inter}.
Then, similarly to \eqref{d-exp},
instead of using Lemma \ref{Lem-U-L2} we apply the Gronwall inequality  to
control the $L^2$-norm of the solution and then obtain the homogeneous local smoothing estimates,
which consequently implies the inhomogeneous estimates \eqref{prop-stri}
by analogous arguments as in the proof of Theorem \ref{Thm-Stri-RNLS}. \hfill $\square$.

We conclude this section with a simplified proof without duality argument for \eqref{Stri-X}
in the conservative case,
but with $(s_1,p_1,q_1) = (s_2,p_2,q_2)$.
In fact, applying Lemma \ref{Lem-stsm-P} and \eqref{bu} to \eqref{equa-rnls2} we get
\begin{align*}
   \|u\|_{\mathcal{X}_{T,s,p,q}}
   \leq& C(1+T)(|X_0|_2 + \|b(\cdot,\cdot, D) u\|_{L^2(0,T; H^{-\frac {m-1}{2}}_{1})}   + \|f\|_{\mathcal{X}'_{T,-s,p',q'}}) \nonumber \\
   \leq& C(1+T)(|X_0|_2 + (1+\beta^*_T+T)^m \| u\|_{L^2(0,T; H^{\frac{m-1}{2}}_{-1})}  + \|f\|_{\mathcal{X}'_{T,-s,p',q'}}),
\end{align*}
where $b(t,x,D) : =  i e^{-W(t,x)}  [\Psi_{P},  e^{W(t,x)} ]$,
$f= e^{-W} F$.
Then, by \eqref{esti-H12-F},
\begin{align*}
   \|u\|_{\mathcal{X}_{T,s,p,q}}
   \leq& C(1+T)^{\frac 32}(1+\beta^*_T + T)^{\frac{(m-2)^2+m}{2}+m} \|u\|_{L^\9(0,T; L^2)}  \nonumber \\
       & + (1+T) (1+\beta^*_T +T)^m  \|f\|_{\mathcal{X}'_{T,-s,p',q'}}.
\end{align*}
Therefore,
similarly to \eqref{v-l2-con},
since $\|u\|_{L^\9(0,T; L^2)} \leq 2 \|u\|^{1/2}_{\mathcal{X}_{T,s,p,q}} \|f\|^{1/2}_{\mathcal{X}'_{T,-s,p',q'}}$,
using Cauchy's inequality $ab\leq \ve a^2 + \ve^{-1} b^2$
we obtain similar estimate as \eqref{Stri-X} but with $(s_1,p_1,q_1) = (s_2,p_2,q_2)$.

\section{Applications} \label{Sec-Appli}

This section contains several applications to nonlinear problems,
including well-posedness, integrability of global solutions
and large deviation principle.

\subsection{Stochastic nonlinear Schr\"odinger equation with variable coefficients} \label{Sec-SNLS}

Consider the stochastic nonlinear Schr\"odinger equation with variable coefficients
and lower order perturbations
\begin{align} \label{equa-qSNLS}
    &dX(t) = i \sum\limits_{j,k=1}^d D_j a^{jk}(x) D_k X(t) dt + b(t,x)\cdot DX(t)  dt + c(t,x)X(t)  dt \nonumber \\
    & \qquad \  - \lbb i |X|^{\a-1}X(t)  dt -\mu X(t)  dt + X (t) dW(t) ,  \nonumber \\
    & X(0)= X_0 \in L^2,
\end{align}
where $a^{jk}$ are real valued and symmetric, $1\leq j,k\leq d$,
$W$ and $\mu$ are as in \eqref{equa-dis},
$b(t,x), c(t,x)$, $t\geq 0$, are continuous $\{\mathcal{F}_t\}$-adapted processes in $\mathbb{C}^d$ and $\mathbb{C}$ respectively,
the coefficient $\lbb =1$ (resp. $\lbb =-1$)
corresponds to the focusing (resp. defocusing) case
and $\a\in (1,\9)$.
We assume that
\begin{enumerate}
  \item[(B1)] {\it Ellipticity}. There exists $C>0$ such that
  \begin{align*}
      C^{-1}|\xi|^2 \leq \<a^{jk}\xi,\xi\> \leq C|\xi|^2.
  \end{align*}
  \item[(B2)] {\it Asymptotic flatness}. For any multi-index $\beta \not =0$,
  \begin{align*}
      | \p_x^\beta a^{jk}(x)| \leq C_\beta \<x\>^{-2}
\end{align*}
  and for $1\leq k,l\leq d$,
  \begin{align}  \label{a-small}
    &\sum\limits_{j\in \mathbb{Z}} \sup\limits_{\{2^j\leq |x|\leq 2^{j+1}\}}
    |x|^2( |\p_{x_k x_l} a(x)| + |x||\p_{x_k} a(x)| + |a(x)-I_n|) \leq \ve \ll 1.
  \end{align}
  Moreover, for any multi-index $\a$,
  \begin{align*}
     \sup\limits_{t\in [0,T]} (| \p_x^\a b(t,x)| +| \p_x^\a c(t,x)|)  \leq g(T) C_\a \<x\>^{-2}, \ \ \bbp-a.s.,
  \end{align*}
  where $g(t)$, $t\geq 0$,  is an $\{\mathcal{F}_t\}$-adapted continuous  process.
\end{enumerate}

\begin{theorem} \label{Thm-QSNLS}
Assume $(B1)$, $(B2)$ and the asymptotic flatness of $e_j$ in \eqref{e-decay}.
Let $\a \in [1,1+4/d]$ and $X_0$ be $\mathcal{F}_0$-measurable,
and $X_0\in L^2$, $\bbp$-a.s..

Then, there exist a stopping time $\tau(\leq T)$ and
a unique solution $X$ to  \eqref{equa-qSNLS} on $[0,\tau]$,
such that
$X\in C([0,\tau]; L^2) \cap \mathcal{X}_{\tau,p,q}$
for any Strichartz pair $(p,q)$, $\bbp$-a.s.
Moreover,  $\tau =T$, $\bbp$-a.s., if $\a \in [1,1+4/d)$
and $b, c$ vanish.
\end{theorem}

{\bf Proof of Theorem \ref{Thm-QSNLS}.}
It follows from \cite{T08} that
the smallness condition \eqref{a-small}
precludes the existence of trapped bicharacteristics and
Assumptions $(A2)$ and $(A3)$  hold for
the operator $P(x,D)=D_j a^{jk}(x) D_k$.
Thus, Theorem \ref{Thm-Stri-per} yields
pathwise Strichartz and local smoothing estimates for
the linear part of equation \eqref{equa-qSNLS}.
Consequently,
similar arguments as in the proof of \cite[Lemma 4.2]{BRZ14}
yield the local well-posedness.
In the mass subcritical case where
$\a \in [1,1+4/d)$  and $b, c$ vanish,
similarly to \eqref{mar-Y},
we have the martingale property of $|X(t)|_2^2$.
Thus, arguing as in the proof of \cite[Proposition 3.2]{BRZ14} we obtain the global existence of solutions
to \eqref{equa-qSNLS}, i.e. $\tau=T$, $\bbp$-a.s.
\hfill $\square$

\subsection{Integrability of global solutions}

Below we consider the typical stochastic nonlinear Schr\"odinger equation
with power nonlinearity
as in \cite{BRZ14}--\cite{BRZ16.2}, namely,
\begin{align} \label{equax}
      &idX =\Delta X dt+\lambda|X|^{\alpha-1}Xdt
      -i\mu Xdt+iXdW,
\end{align}
with $X(0)=X_0$ being $\mathcal{F}_0$-measurable. As mentioned in Section \ref{Sec-Intro},
global well-posedness of \eqref{equax} was first  studied in \cite{BD99,BD03}.
Pathwise global well-posedness
in the full mass  and energy subcritical cases has been recently obtained in \cite{BRZ14, BRZ16}.
See also \cite{H16} for the full mass subcritical case.

Motivated by optimal control problems (see \cite{BRZ16.2}),
we shall prove  $L^\rho(\Omega)$-integrability of global solutions
in both mass and energy subcritical cases,
which can be viewed as a complement to \cite{BRZ14,BRZ16}.

\begin{theorem} \label{Thm-SNLS}
Assume the asymptotically flat condition of $e_j$, $1\leq j\leq N$, as in \eqref{e-decay}.
Let $\a\in (1,1+4/d)$ (resp. $\a\in (1,1+4/d)$ if $\lbb =1$,
and $\a\in (1, 1+ 4/(d-2)_+)$ if $\lbb =-1$)
and $X_0\in L^\rho(\Omega; L^2)$ (resp. $L^\rho(\Omega; H^1)$) for any $1\leq \rho <\9$.

Then, for each $0<T<\9$,
there exists a unique global   solution $X$  to \eqref{equax} on $[0,T]$,
such that $X\in C([0,T]; L^2)$ (resp. $X\in C([0,T]; H^1)$), $\bbp$-a.s.,
and  for any Strichartz pair $(p,q)$ and $1\leq \rho <\9$,
\begin{align} \label{Integ-Xpq}
      \bbe \|X\|_{\calx_{T,p,q}}^\rho  <\9 \ \
     (resp.\  \bbe  (\|X\|^\rho_{\calx_{T,p,q}} + \|\na X\|^\rho_{\calx_{T,p,q}})  <\9 ).
\end{align}
\end{theorem}

{\bf Proof.}
The global well-posedness follows from
similar arguments as in \cite{BRZ14,BRZ16}.
For the integrability of solutions,
let us first consider the $L^2$ case.

Choose the Strichartz pair $(p,q) = (\a+1, \frac{4(\a+1)}{d(\a-1)})$
and set $u= e^{-\Phi(W)}X$.
As in the proof of \cite[(4.9), (4.10)]{BRZ14},
applying Theorem \ref{Thm-Stri-L2} and  H\"older's inequality we have
that for any Strichartz pair $(p,q)$,
\begin{align} \label{Integ-Xpq*}
   \|u\|_{\calx_{t,p,q}}
   \leq C_T (|X_0|_2 +  t^\theta \g_T \|u\|^\a_{L^q(0,t; L^p)}),
\end{align}
where $\g_T = e^{(\a-1)\|\Phi(W)\|_{L^\9(0,T; L^\9)}} $,
$\theta = 1-  d(\a-1)/4 \in (0,1)$,
and $C_T \in L^\rho (\Omega)$ for any $1\leq \rho <\9$.
Then, taking
\begin{align*}
    t =(\a-1)^{\frac{\a-1}{\theta}} \a^{-\frac{\a}{\theta}}(|X_0|_2+1)^{-\frac{\a-1}{\theta}}  C_T^{-\frac{\a}{\theta}} \g_T^{-\frac 1 \theta} (\leq T),
\end{align*}
and using \cite[Lemma 6.1]{BRZ16.2} we get
\begin{align*}
    \|u\|_{\calx_{t,p,q}}
    \leq \frac{\a}{\a-1} C_T |X_0|_2 \leq\frac{\a}{\a-1}  C_T \|u\|_{C([0,T]; L^2)}.
\end{align*}
Iterating similar arguments on $[jt, (j+1)t\wedge T]$, $1\leq j\leq [T/t]$,
since $1/q < 1/2 $, we obtain
\begin{align} \label{yT}
     \|u\|_{\calx_{T,p,q}}
     \leq& \frac{2\a}{\a-1} C_T ([\frac T t]+1)^{\frac 12}\|u\|_{C([0,T]; L^2)}  \nonumber \\
     \leq& \frac{2\a}{\a-1} C_T ([\frac T t]+1)^{\frac 12}\|e^{-\Phi(W)}\|_{L^\9(0,T; L^\9)}  \|X\|_{C([0,T]; L^2)}.
\end{align}
Therefore, since
$\|X\|_{\calx_{T,p,q}}
   \leq C(l) \sup_{t\in [0,T]} |e^{\Phi(W)(t)}|_{S^0}^{(l)} \|u\|_{\calx_{T,p,q}} $
for some $l \geq 1$,
using the $L^\rho(\Omega)$-integrability of $C_T$,
$e^{\|\Phi(W)\|_{L^\9(0,T;L^\9)}}$ and
$\|X\|_{C([0,T]; L^2)}$
(see \cite[Lemma $3.6$]{BRZ16})
we obtain \eqref{Integ-Xpq} in the $L^2$ case.

The proof in the $H^1$ case is similar. Indeed,
as in the proof
of \cite[(2.25)]{BRZ16}, using Theorem \ref{Thm-Stri-L2}
and the Sobolev imbedding $|u|_{L^p} \leq D|u|_{H^1}$ we have
\begin{align*}
   \|u\|_{\calx_{t, p,q}}+ \|\na u\|_{\calx_{t, p,q}}
   \leq C_T (|X_0|_{H^1} +  t^\theta D(T) \|u\|^{\a-1}_{C([0,T]; H^1)} \|u\|_{L^q(0,t; W^{1,p})}),
\end{align*}
where $\theta = 1 - 2/q \in(0,1)$,
and $D(T) = \a D^{\a-1} (\|\na \Phi(W)\|_{L^\9(0,T; L^\9)}+ 2) \g_T$.
Taking $t= (2C_TD(T)\|y\|^{\a-1}_{C([0,T]; H^1)})^{-1/ \theta}$
and using iterated arguments we get
\begin{align} \label{naT}
   \|u\|_{\calx_{T,p,q}} + \|\na u\|_{\calx_{T,p,q}} \leq 8   C_T ([\frac T t]+1)^{\frac 12} \|u\|_{C([0,T];H^1)}.
\end{align}
Since $\|X\|_{C([0,T]; H^1)}$ is $L^\rho(\Omega)$-integrable (see \cite[(2.3)]{BRZ16.2}),
and so are the coefficients $C_T$ and $\g_T$,
we obtain \eqref{Integ-Xpq} in the $H^1$ case. \hfill $\square$

\subsection{Large deviation principle}

We first consider the large deviation principle (LDP) for the small noise asymptotics
for \eqref{equa-dis} in the conservative case.
Consider
\begin{align} \label{equa-disve}
      dX^\ve(t)= i P(x,D) X^\ve(t)dt  + F(t) dt  - \ve \mu X^\ve(t)dt + \sqrt{\ve} X^\ve(t)dW,
\end{align}
where $X^\ve(0)=X_0 \in L^2$, $\bbp$-a.s., $X_0$ is $\mathcal{F}_0$-measurable,
$\mu$ and $W$ are as in \eqref{W},
${\rm Re} \mu_j =0$, $1\leq j\leq N$, $\ve\in (0,1)$,  and
$F$ is $\{\mathcal{F}_t\}$-adapted, $F\in \mathcal{X}'_{T,-s,p',q'}$ for some $(s,p,q)\in \mathcal{A}$, $\bbp$-a.s.

Let $C_0([0,T]; \bbr^N) =\{u\in C([0,T]; \bbr^N): u(0)=0\}$.
Introduce
the map $\mathcal{G}: C_0([0,T];\bbr^N) \to \mathcal{X}_{T,s,p,q}$,
such that for any $g = (g_1,\cdots,g_N)  \in C_0([0,T];\bbr^N)$,
$u^g := \mathcal{G}(g)$ solves the equation
\begin{align} \label{equa-ug}
     \partial_t u^g (t) = i e^{-\wt{\Phi}(g)} P(x,D) e^{\wt{\Phi}(g)} u^g(t)  + e^{-\wt{\Phi}(g)} F(t),
\end{align}
where $u^g(0) = X_0$,
and $\wt{\Phi}(g) =  \sum_{j=1}^N \mu_j e_j g_j$.

Moreover, define the map
$\mathcal{S}: C_0([0,T];\bbr^N) \to \mathcal{X}_{T,s,p,q}$ by
\begin{align} \label{S-G}
   \mathcal{S}(g) = e^{\wt{\Phi}(g)} \mathcal{G}(g), \ \ \forall g \in C_0([0,T]; \bbr^N).
\end{align}

Thus, by the rescaling \eqref{rescal},
\begin{align} \label{X-beta}
    X^\ve = \mathcal{S}(\sqrt{\ve}\beta) = e^{\wt{\Phi}(\sqrt{\ve}\beta)} \mathcal{G}(\sqrt{\ve}\beta),
\end{align}
where $\beta = (\beta_1,\cdots, \beta_N)$ are $N$ dimensional real valued Brownian motions.

\begin{theorem} \label{Thm-LDP}
The family $\{X^\ve\}$ satisfies a LDP on $\calx_{T,s,p,q}$ of speed $\ve$ and a good rate function
\begin{align} \label{LDP-dis}
   I (w) = \frac 12 \inf\limits_{g\in H^1(0,T; \bbr^N): w= \mathcal{S}(g)} \|\dot{g}\|^2_{L^2(0,T; \bbr^N)},
\end{align}
where $\dot{g}$ denotes the derivative of $g$.
\end{theorem}

The key observation here is, that
the solution map $\mathcal{G}$
of the reduced equation \eqref{equa-ug} is continuous from $ C([0,T]; \bbr^N)$ to $\mathcal{X}_{T,s,p,q}$,
i.e., the solution to \eqref{equa-ug} depends continuously on lower order perturbations.
This fact implies, via the representation formula \eqref{X-beta} of the stochastic solution to \eqref{equa-disve},
the large deviation principle for  $\mathcal{S}$
by virtue of Varadhan's contraction principle.

\begin{lemma} \label{Lem-yn}
The map $\mathcal{G}: C_0([0,T]; \bbr^N) \mapsto \mathcal{X}_{T,s,p,q}$ is continuous.
\end{lemma}

{\bf Proof.}
Let $g_n, g\in C_0([0,T]; \bbr^N)$, $g_n \to g$ in $ C_0([0,T]; \bbr^N)$, as $n\to \9$.
Set $u_n = \mathcal{G}(g_n)$.
Then,  by \eqref{equa-ug},
if $b(x,D,g_n) := i e^{-\wt{\Phi}(g_n)} [P(x,D), e^{\wt{\Phi}(g_n)}]$,
\begin{align*}
  \partial_t u_n (t) = i P(x,D) u_n  + b(x,D,g_n) u_n + e^{-\wt{\Phi}(g_n)} F(t).
\end{align*}
Define $u,b(x,D,g)$ similarly as above. Then,
\begin{align} \label{LDP-dis*}
   \p_t(u_n-u)
   =& i P(u_n-u) + (b(x,D,g_n)u_n - b(x,D,g)u) + (e^{-\wt{\Phi}(g_n)} -e^{-\wt{\Phi}(g)})F(t) \nonumber \\
   =& i e^{-\wt{\Phi}(g_n)} P(x,D) e^{\wt{\Phi}(g_n)} (u_n-u) + (b(x,D,g_n)-b(x,D,g)) u \nonumber \\
    & + (e^{-\wt{\Phi}(g_n)} -e^{-\wt{\Phi}(g)})F(t).
\end{align}

Note that, Strichartz and local smoothing estimates as in \eqref{prop-stri} also holds for
$e^{-\wt{\Phi}(g_n)} P(x,D) e^{\wt{\Phi}(g_n)}$
with $e^{-\Phi(W)}X$ and $e^{-\Phi(W)} F$ replaced by
$u_n$ and $e^{-\wt{\Phi}(g_n)} F$ respectively,
and the constant $C_T$ is independent of $n$,
due to the boundedness of $\sup_n\|g_n(t)\|_{C([0,T];\bbr^N)}$.
Then,
taking into account $\<x\>\<D\>^{-(m-1)/2}(b(x,D,g_n)
-b(x,D,g))$ $\<D\>^{-(m-1)/2}\<x\> \in S^0$,
and using Lemma \ref{Lem-L2-Bdd}
we obtain
\begin{align*}
       & \|u_n-u\|_{\mathcal{X}_{T,s,p,q}} \\
   \leq& C_T \|(b(x,D,g_n)-b(x,D,g))u\|_{L^2(0,T; H^{-\frac{m-1}{2}}_{1})}
         + C_T \|(e^{-\wt{\Phi}(g_n)} -e^{-\wt{\Phi}(g)})F\|_{\mathcal{X}'_{T,-s,p',q'}}  \\
   \leq& C_T \sup\limits_{t\in[0,T]} |b(x,\xi,g_n(t))-b(x,\xi,g(t))|_{S^{m-1}}^{(l)} \|u\|_{L^2(0,T; H^{\frac{m-1}{2}}_{-1})} \\
       & + C_T \sup\limits_{t\in[0,T]} |e^{-\wt{\Phi}(g_n)(t)} -e^{-\wt{\Phi}(g)(t)}|_{S^0}^{(l)} \|F\|_{\mathcal{X}'_{T,-s,p',q'}}
\end{align*}
for some $l\geq 1$.

Thus, since
\begin{align*}
   b(x,\xi,g_n)-b(x,\xi,g)
   = i\sum\limits_{1\leq |\a|\leq m} \frac{1}{\a !} \p^\a_\xi P ( \psi_\a(\wt{\Phi}(g_n)) - \psi_\a(\wt{\Phi}(g))),
\end{align*}
where $ \psi_\a(\wt{\Phi}(g_n)) = e^{-\wt{\Phi}(g_n)} D_x^\a e^{\wt{\Phi}(g_n)}$
and $\psi_\a(\wt{\Phi}(g))$ is defined similarly,
using the convergence of $\{g_n\}$ we get that
$\sup_{t\in[0,T]}|b(x,\xi,g_n(t))-b(x,\xi,g(t))|_{S^{m-1}}^{(l)} \to 0$,
and  $\sup_{t\in[0,T]} |(e^{-\wt{\Phi}(g_n)(t)} -e^{-\wt{\Phi}(g)(t)})|_{S^0}^{(l)} \to 0$, as $n\to \9$,
which implies that $ \|u_n-u\|_{\mathcal{X}_{T,s,p,q}} \to 0$,
thereby completing the proof. \hfill $\square$

\begin{corollary} \label{Cor-Xn}
The map $\mathcal{S}: C([0,T];\bbr^N) \to \mathcal{X}_{T,s,p,q}$ is continuous.
\end{corollary}

{\bf Proof.}
Let $g_n, g, u_n, u$ be as in the proof of Lemma \ref{Lem-yn},
and set  $X_n= \mathcal{S}(g_n)$, $X= \mathcal{S}(g)$.
Then, by \eqref{S-G},
\begin{align*}
   \|X_n - X\|_{\mathcal{X}_{T,s,p,q}}
   \leq \| e^{\wt{\Phi}(g_n)} (u_n-u) \|_{\mathcal{X}_{T,s,p,q}}
        +  \|(e^{\wt{\Phi}(g_n)} - e^{\wt{\Phi}(g)})u \|_{\mathcal{X}_{T,s,p,q}} .
\end{align*}
Similarly to Lemma \ref{Lem-Psi-Esti},
$ \|e^{\wt{\Phi}(g_n)} \|_{\mathcal{L}(\mathcal{X}_{T,s,p,q})} \leq C |e^{\wt{\Phi}(g_n)}|_{S^0}^{(l)}$
for some $l\in \mathbb{N}$.
Thus, using Lemma \ref{Lem-yn} and the convergence of $g_n$ we obtain
\begin{align*}
   \|X_n - X\|_{\mathcal{X}_{T,s,p,q}}
   \leq& C \sup\limits_{n\in\mathbb{N}}   \sup\limits_{t\in [0,T]}|e^{\wt{\Phi}(g_n)(t)}|_{S^0}^{(l)} \| u_n-u \|_{\mathcal{X}_{T,s,p,q}}  \\
       & + C \sup\limits_{t\in [0,T]} |e^{\wt{\Phi}(g_n)(t)} - e^{\wt{\Phi}(g)(t)}|_{S^0}^{(l)} \|u\| _{\mathcal{X}_{T,s,p,q}} \to 0,
\end{align*}
which finishes the proof. \hfill $\square$

{\bf Proof of Theorem \ref{Thm-LDP}.}
By Schilder's theorem (see e.g. \cite[Theorem 5.2.3]{DZ10}),
$\{\sqrt{\ve} \beta \}$ satisfies the LDP of speed $\ve$
and the good rate function
$$I^\beta = \frac 12 \inf\limits_{g\in H^1(0,T; \bbr^N)} \|\dot{g}\|^2_{L^2(0,T; \bbr^N)}. $$
Then, by virtue of the continuity of $\mathcal{S}$ in Corollary \ref{Cor-Xn}
and Varadhan's  contraction principle (see \cite[Theorem 4.2.1]{DZ10}) we prove  \eqref{LDP-dis}.
\hfill $\square$\\

We conclude this section with the large deviation principle
for the nonlinear
Schr\"odinger equation \eqref{equa-qSNLS}
with variable coefficients in the case where $b$ and $c$ vanish.
See also \cite{G05} for the case of constant coefficients.

As above, for any Strichartz pair $(p,q)$,
introduce the map $\wt{\mathcal{G}}: C_0([0,T]; \bbr^N)$ $\to \mathcal{X}_{T,p,q}$,
such that for any $g\in C_0([0,T];\bbr^N)$, $u^g= \wt{\cal{G}} (g)$ solves the equation
\begin{align} \label{equa-u-SNLS}
     \partial_t u^g  = i e^{-\wt{\Phi}(g)} P(x,D) e^{\wt{\Phi}(g)} u^g - \lbb i|u^g|^{\a-1} u^g,
\end{align}
and $u^g(0)=X_0$,
where $P(x,D) = \sum_{j,k=1}^d D_j a^{jk}(x) D_k$.

Set $\wt{\mathcal{S}}(g) := e^{\wt{\Phi}(g)} \wt{\mathcal{G}}(g)$, $g\in C_0([0,T]; \bbr^N)$.
Then,  $X^\ve = \wt{\mathcal{S}}(\sqrt{\ve} \beta)$ solves \eqref{equa-qSNLS}
with $b,c=0$ and $W$, $\mu$ replaced by $\sqrt{\ve} W$
and $\ve \mu$ respectively.

\begin{theorem} \label{Thm-LDP-SNLS}
Assume the conditions of Theorem \ref{Thm-QSNLS} to hold.
Assume in addition that $\a\in (1,1+4/d)$, $b$ and $c$ vanish, and ${\rm Re} \mu_j =0$, $1\leq j\leq N$.
Then, for any Strichartz pair $(p,q)$,
the family $\{X^\ve\}$ satisfies the LDP on $\mathcal{X}_{T,p,q}$
of speed $\ve$ and a good rate function
\begin{align*}
   \wt{I}(w) = \frac 12 \inf\limits_{g\in H^1(0,T; \bbr^N), w=\wt{S}(g)} \|\dot{g}\|^2_{L^2(0,T; \bbr^N)},
\end{align*}
where $\dot{g}$ denotes the derivative of $g$.
\end{theorem}

\begin{lemma}
Assume the conditions of Theorem \ref{Thm-LDP-SNLS} to hold.
Then,
for any Strichartz pair $(p,q)$,
the map $\wt{\mathcal{G}} : C_0([0,T]; \bbr^N) \to \mathcal{X}_{T,p,q}$ is continuous.
\end{lemma}

{\bf Proof.} Let $g_n, g$ be as in the proof of Lemma \ref{Lem-yn}.
Set $u_n = \wt{\mathcal{G}} (g_n)$, $u= \wt{\mathcal{G}}(g)$ and
choose the Strichartz pair $(p_0,q_0) = (\a+1, \frac{4(\a+1)}{d(\a-1)})$.

We first claim that,
\begin{align} \label{Bdd-un}
     \sup\limits_{n\geq 1} \|u_n\|_{\mathcal{X}_{T,p_0,q_0}} < \9.
\end{align}
In fact,
similarly to \eqref{Integ-Xpq*},
\begin{align*}
   \|u_n\|_{\calx_{t,p_0,q_0}}
   \leq C_T (|X_0|_2 + t^\theta \|u_n\|^\a_{L^{q_0}(0,t; L^{p_0})}),\ \ \forall t\in (0,T),
\end{align*}
where $\theta = 1-d(\a-1)/4 \in (0,1)$,
$C_T$ is independent of $n$.
As in the proof of \cite[(3.16)]{BRZ16.2}, taking
$t= \a^{-\a/\theta} (\a-1)^{(\a-1)/\theta} (|X_0|_2 +1)^{-(\a-1)/\theta} C_T^{-\a/\theta} (\leq T) $,
we have
\begin{align*}
     \|u_n\|_{\calx_{t,p_0,q_0}}  \leq \frac{\a}{\a-1} C_T |X_0|_2.
\end{align*}
Since $|u_n(t)|_2 = |X_0|_2$, $\forall t\in [0,T]$,
we can iterate similar arguments on $[jt, (j+1)t\wedge T]$, $1\leq j\leq [T/t]$, and obtain
\begin{align*}
    \|u_n\|_{\calx_{T,p_0,q_0}} \leq  \frac{2\a}{\a-1} C_T ([\frac T t]+1)^{\frac 1q} |X_0|_2,\ \ \forall n\geq 1,
\end{align*}
which implies \eqref{Bdd-un}, as claimed.

Now, note that, similarly to \eqref{LDP-dis*},
\begin{align*}
   \partial_t (u_n -u)
    =&  i e^{-\wt{\Phi}(g_n)} P(x,D) e^{\wt{\Phi}(g_n)} (u_n -u)
        + (b(x,D,g_n)- b(x,D,g)) u  \\
     &   - \lbb i (|u_n|^{\a-1} u_n - |u|^{\a-1}u),
\end{align*}
where $b(x,D,g_n)$ and $b(x,D,g)$ are defined as in \eqref{LDP-dis*}.

Applying Strichartz estimates and H\"older's inequality we obtain
\begin{align*}
    \|u_n - u\|_{\mathcal{X}_{t, p, q}}
    \leq C_T ( w_n(u) + \wt{C} t^\theta \|u_n - u\|_{L^{q_0}(0,t; L^{p_0})} ),
\end{align*}
where $\wt{C} = \a( \sup_{n\geq 1} \|u_n\|^{\a-1}_{L^{q_0}(0,T; L^{p_0})} + \|u\|^{\a-1}_{L^{q_0}(0,T; L^{p_0})}) <\9$,
and
\begin{align*}
   w_n(u) =& \|(b(x,D,g_n)- b(x,D,g)) u\|_{L^2(0,T; H^{-\frac 12}_1)}  \\
          \leq& C \sup\limits_{t\in[0,T]} |b(x,\xi,g_n(t))- b(x,\xi,g(t))|_{S^1}^{(l)} \|u\|_{L^2(0,T; H^{\frac 12}_{-1})}
          \to 0.
\end{align*}
Thus, taking $t= (2C_T \wt{C})^{-1/\theta}$ we obtain
$\|u_n - u\|_{\mathcal{X}_{t, p, q}} \leq 2 C_T w_n(u) \to 0. $

Thus, as $t$ is independent of $n$, iterating similar estimates on $[jt, (j+1)t\wedge T]$,
$1\leq j\leq [T/t]$,
we obtain $ \|u_n - u\|_{\mathcal{X}_{T, p, q}} \to 0$
and complete the proof. \hfill $\square$

Now, using similar arguments as in the proof of Theorem \ref{Thm-LDP} we prove Theorem \ref{Thm-LDP-SNLS}.

\section{Appendix} \label{Sec-App}

{\bf Proof of Lemma \ref{Lem-EstiErr}.} First note that, for any $l,k\geq 1$,
\begin{align} \label{ctheta}
     |c_\theta(x,\xi)|
   =& \bigg|\iint e^{-iz\cdot \eta} \<\eta\>^{-k} \<\p_z\>^k \<z\>^{-l} \<\p_\eta\>^l
      (a(x,\xi+\theta\eta) b(x+z,\xi)) d\eta dz \bigg| \nonumber \\
   \leq& C C_1(l)C_2(k) \<x\>^{\rho_1(0)} \<\xi\>^{m_2}
        \int  K_1(\xi, \eta)  d\eta
        \int  K_2(x,z) dz,
\end{align}
where $K_1(\xi, \eta)=\<\eta\>^{-k} \<\xi+\theta \eta\>^{m_1-l} $ and
$K_2(x,z) =  \<z\>^{-l} \<x+z\>^{\rho_2(0)}$.

We first consider $K_1(\xi,\eta)$.
Note that,
for $\eta\in \Lambda : =\{\eta: |\eta|<\<\xi\>/2\}$,
$\<\xi\>/2 \leq  \<\xi+\theta \eta\> \leq 3\<\xi\>/2$. Then,
\begin{align*}
    \int\limits_{\Lambda} K_1(\xi,\eta) d\eta
    \leq C \int\limits_{\Lambda} \<\eta\>^{-k} \<\xi\>^{m_1-l} d\eta
    \leq C  \<\xi\>^{m_1-l+d},
\end{align*}
which implies that for $l>d$, $\int_{\Lambda} K_1(\xi,\eta) d\eta  \leq C \<\xi\>^{m_1}$.
Moreover, for $\eta \in  \Lambda^c$,
$\<\xi+\theta \eta\> \leq \<\xi\>+|\eta| \leq 3 |\eta|$.
Since $\<\eta\> \geq  |\eta|$, we obtain
\begin{align*}
    \int\limits_{\Lambda^c} K_1(\xi,\eta) d\eta
    \leq C \int\limits_{\Lambda^c} |\eta|^{-k+(m_1)_+} d\eta
    \leq C \<\xi\>^{-k+(m_1)_+ +d}
    \leq C \<\xi\>^{m_1},
\end{align*}
where we choosed $k$ such that $-k+(m_1)_+ +d < - (m_1)_-$,
where $(m_1)_+ = \max\{m_1, 0\}$, $(m_1)_- = (-m_1)_+$.
Thus, for $l>d$, $k>|m_1| +d $, we have
\begin{align} \label{K1}
     \int K_1(\xi, \eta) d\eta \leq C \<\xi\>^{m_1}.
\end{align}

The estimate for $K_2(x,z)$ is similar.
Set $\Omega :=\{z:|z|\leq \<x\>/2 \}$.
For $z\in \Omega$, $\<x\> /2 \leq \<x+z\> \leq 3\<x\>/2 $, and so
\begin{align*}
  \int \limits_{\Omega} K_2(x,z) dz
  \leq C \<x\>^{\rho_2(0)} \int \limits_{\Omega}\<z\>^{-l} dz
  \leq C \<x\>^{\rho_2(0)},
\end{align*}
if $l>d$. Moreover, for $z  \in \Omega^c$,
$\<x+z\> \leq 3|z|$. Then,
\begin{align*}
   \int \limits_{\Omega^c} K_2(x,z) dz
   \leq C \int\limits_{\Omega^c} |z|^{-l+(\rho_2(0))_+} dz
   \leq C \<x\>^{-l+(\rho_2(0))_+ +d}
   \leq C\<x\>^{\rho_2(0)},
\end{align*}
if $l$ is large enough such that $-l+(\rho_2(0))_+ + d < -(\rho_2(0))_- $.
Thus,
for $l>|\rho_2(0)|+d$, we have
\begin{align} \label{K2}
    \int K_2(x,z)dz \leq C\<x\>^{\rho_2(0)}.
\end{align}

Therefore, plugging \eqref{K1} and \eqref{K2} into \eqref{ctheta} we obtain \eqref{esti-ctheta}. \hfill $\square$

{\it \bf Proof of Lemma \ref{Lem-Psi-Esti}.}
For simplicity, we set $a_R (x,\xi) := a(x,\xi) \theta_R(\xi)$.

$(i)$. By straightforward computations, for any $l\in \mathbb{N}$,
\begin{align} \label{e-s0l}
     |c_R|_{S^0}^{(l)} +  |c_R^{-1}|_{S^0}^{(l)}
     \leq C(l)  M^l e^{M \|a\|_{\9}},
\end{align}
which along with Lemma \ref{Lem-L2-Bdd} implies \eqref{esti-Psi-1-1}.

$(ii)$. By Lemma \ref{Lem-Err}, $\Psi_{c_R} \Psi_{c^{-1}_R} = I + \Psi_{e_R}$ with
\begin{align*}
   e_R = \frac{1}{(2\pi)^d}\int_0^1 \sum\limits_{|\g|=1}  \iint
          e^{iy\cdot \eta} \p_\xi^\g c_R(x,\xi+\theta \eta) D_x^\g c_R^{-1}(x+y,\xi) dy d\eta d\theta \in S^0.
\end{align*}
Using Corollary \eqref{Cor-EstiErr} we have that
for any $l\in \mathbb{N}$,
$|e_R|_{S^0}^{(l)}
    \leq C(l') \sum_{|\a|=1} |\partial^\a_\xi c_R|_{S^0}^{(l')}$ $|\partial_x^\a c_R^{-1}|_{S^0}^{(l')}$
for some $l'\in\mathbb{N}$.
Below we shall prove that
\begin{align} \label{crcr-1}
    |\partial^\a_\xi c_R|_{S^0}^{(l')}  |\partial_x^\a c^{-1}_R|_{S^0}^{(l')}
    \leq C(l') R^{-1} M^{2l'} e^{2M\|a\|_{\9}}
\end{align}
Then, in view of Lemma \ref{Lem-L2-Bdd},
for $R= C(l) M^{2l}e^{2M\|a\|_{\9}}$ with $l$ and $C(l)$ large enough,
$\|\Psi_{e_R}\|_{\mathcal{L}^2(L^2)}$ is less than $1/2$,
which yields that  $I+ \Psi_{e_R}$ is invertible and $\Psi_{c_R}^{-1} = \Psi_{c^{-1}_R} (I+\Psi_{e_R})^{-1}$,
thereby implying \eqref{esti-Psi-inver} by  \eqref{e-s0l}.

It remains to prove \eqref{crcr-1}. Note that
\begin{align}  \label{esti-pcR}
     | \partial_\xi^{\a} c_R |^{(l)}_{S^0}
     = M | c_R \partial^\a_\xi a_R|^{(l)}_{S^0}
     \leq  M |c_R|^{(l)}_{S^0} |\partial^\a_\xi a_R|^{(l)}_{S^0}
     \leq  C(l) M R^{-1}|c_R|^{(l)}_{S^0}.
\end{align}
Similarly, since $|\partial^\a_x a_R|^{(l)}_{S^0}
    = | \theta_R \partial^\a_x a|^{(l)}_{S^0}
    \leq C(l)$,  we have
\begin{align} \label{pcr}
     | \partial_x^{\a} c^{-1}_R |^{(l)}_{S^0}
     = |M c_R^{-1} \partial_x^\a a_R|^{(l)}_{S^0}
     \leq& C(l) M  |c^{-1}_R|^{(l)}_{S^0}.
\end{align}
Thus, \eqref{crcr-1} follows from \eqref{e-s0l}, \eqref{esti-pcR} and \eqref{pcr},
and \eqref{esti-Psi-inver} is proved.

$(iii)$ First note that  for any $f\in H^{\frac{m-1}{2}}_{-1}$,
\begin{align*}
   \|\Psi_{c_R} f\|_{H^{\frac{m-1}{2}}_{-1}}
  =  |b(x,D)\<x\>^{-1}\<D\>^{\frac{m-1}{2}}f |_2
  \leq   C |c_R|^{(l)}_{S^0} \|f\|_{H^{\frac{m-1}{2}}_{-1}}
\end{align*}
for some $l\in \mathbb{N}$, where
$b(x,D) := \<x\>^{-1}\<D\>^{(m-1)/2} \Psi_{c_R} \<D\>^{-(m-1)/2} \<x\>  \in S^0$
due to  Corollary \ref{Cor-EstiErr}, and we used Lemma \ref{Lem-L2-Bdd}.

Similarly, since
$\Psi^{-1}_{c_R} = \Psi_{c^{-1}_R} (I+\Psi_{e_R})^{-1}$, we have
that
\begin{align*}
   \|\Psi^{-1}_{c_R} f\|_{H^{\frac{m-1}{2}}}
   \leq   C |c^{-1}_R|^{(l_0)}_{S^0} |r_R|^{(l_0)}_{S^0} \|f\|_{H^{\frac{m-1}{2}}_{-1}}, \ \ \forall f\in H^{\frac{m-1}{2}}_{-1},
\end{align*}
for some $l_0\geq 1$, where $r_R\in S^0$ is the symbol of $(I+\Psi_{e_R})^{-1}$.

Below we claim that, for $R= C(l) M^l e^{2M\|a\|_{\9}}$ with $C(l)$ and $l$ large enough,
there exists $C$ independent of $M$ and $R$, such that
\begin{align} \label{rR}
    |r_R|^{(l_0)}_{S^0} \leq C.
\end{align}
For this purpose, let $e_k\in S^0$ be the symbol of $\Psi^k_{e_R}$.
From the proof of $(ii)$ we see that
\begin{align*}
    \Psi_{r_R} = (I+\Psi_{e_R})^{-1}
    = \sum\limits_{k=0}^\9 (-1)^k \Psi^k_{e_R}
    = \sum\limits_{k=0}^\9 (-1)^k \Psi_{e_k},
\end{align*}
which implies that $r_R = \sum_{k=0}^\9 (-1)^k e_k$.
Note that, for any multi-indices $\a,\beta$, $|\a+\beta| \leq l_0$,
by Corollary \ref{Cor-EstiErr}, \eqref{esti-pcR} and \eqref{pcr},
there exits $l\geq 1$ such that
\begin{align*}
   |\p_\xi^\a \p_x^\beta e_R(x,\xi) |
\leq C(l) |c_R|^{(l)}_{S^0} |c^{-1}_R|^{(l)}_{S^0} M^2 R^{-1}  \<\xi\>^{-|\a|}
    = \ve(R,l)   \<\xi\>^{-|\a|},
\end{align*}
where $\ve(R,l) =  C(l) |c_R|^{(l)}_{S^0} |c^{-1}_R|^{(l)}_{S^0} M^2 R^{-1} \to 0$, as $R\to \9$.
Applying Corollary \ref{Cor-EstiErr} again we have
$|\p_\xi^\a \p_x^\beta e_k(x,\xi)|
    \leq C^k(l_0)  \ve^k(R,l)  \<\xi\>^{-|\a|}$.

Therefore, for $R\geq 2 C(l_0) C(l) |c_R|^{(l)}_{S^0} |c^{-1}_R|^{(l)}_{S^0} M^2$,
we obtain
\begin{align*}
     |\p_\xi^\a \p_x^\beta r_R(x,\xi) |
   = \bigg| \sum\limits_{k=0}^\9 (-1)^k \p_\xi^\a \p_x^\beta e_k(x,\xi) \bigg|
  \leq 2 \<\xi\>^{-|\a|},
\end{align*}
which implies \eqref{rR} as claimed and so proves  \eqref{esti-Psi-h12}.
The proof is complete. \hfill $\square$ \\

{\it \bf Acknowledgements.}
The author is grateful to Professor Michael R\"ockner for helpful discussions.
This work is supported by NSFC (No. 11501362).
Financial support through SFB 1283 at Bielefeld University is also gratefully acknowledged.

\end{document}